\documentclass[12pt]{article}

\newtheorem{lemma}{Lemma}
\newtheorem{theorem}{Theorem}

\newtheorem{definition}{Definition}
\newtheorem{remark}{Remark}

\usepackage{amsmath,amssymb,ascmac,color,ulem}

\usepackage{amsmath,amssymb,ascmac,color,ulem}

\begin{document}

\title{ Simple Singularity of type $E_7$\\
 and\\
  the Complex Reflection Group ST34}
\author{Jiro Sekiguchi\thanks{Department of Mathematics, Tokyo University of Agriculture and Technology}}

\maketitle

\begin{flushright}
\end{flushright}
\begin{abstract}
This paper studies a family of surfaces of ${\bf C}^3$ which is a  deformation of a
simple singularity of type $E_7$.
This family has six parameters which are regarded as basic invariants of the complex 
reflection group No.34 in the list of the paper of Shephard and Todd \cite{ST}.
We compute 1-parameter subfamilies of the family in question
corresponding to
 corank one reflection subgroups of No.34 group.
 In particular, we  determine the 
types of simple singularities on the surfaces appeared in this manner.
\end{abstract}


\section{Introduction}
\label{section:intro}

By the classification of unitary reflection groups of Shephard-Todd \cite{ST},
it is known that there are
 three infinite series of unitary reflection groups, plus thirty-four exceptional  ones  numbered as $4,5,\dots,37$.
The groups numbered  as $35,\>36,\>37$ are real reflection groups of type $E_6,\>E_7,\>E_8$, respectively.
We focus our attention on the group numbered as 34  which we denote ST34 in this paper.
 (ST34 is also called the Mitchell group.)
The order of ST34 is largest in the exceptional complex reflection groups
and among other exceptional groups, ST34 is most difficult to treat.
For example, basic invariants of ST34 are extremely complicated to write down 
compared with other cases (cf. \cite{CS}, \cite{OS}).

In this paper, we shall treat Arnold's problem for the complex reflection group ST34 and related topics.
First of all, we recall the formulation of the  problem.
In the book \cite{Ar}, p.20, it is written that

\vspace{5mm}
``1974-5 Find applications of the (Shephard-Todd) complex reflection groups to singularity theory.''

\vspace{5mm}
There are many studies concerning this problem
(cf. comments by V.V. Goryunov in \cite{Ar}, p.305, \cite{Go}, P. Slodowy \cite{Sd2}, T. Yano \cite{Y1}, \cite{Y2},
\cite{SeV}).

In the course of the study of algebraic potentials,
the author recognized the existence of  a family ${\cal F}$ of surfaces of ${\bf C}^3$
which is a  deformation of a simple singularity of type $E_7$ and is also closely
related with the group ST34.
Each member of ${\cal F}$ is a surface in ${\bf C}^3$ with the coordinate $(x,y,z)$ defined by
the polynomial of $(x,y,z)$ depending on $\tau=(t_1,t_2,t_3,t_4,t_5,t_7)\in{\bf C}^6$:
\begin{equation}
\label{equation:int-f}
f_{\tau}(x,y,z)=
y^3+x^3y+t_7x+t_5x^2+t_3x^3+t_1x^4+y(t_4x+t_2x^2)-z^2=0.
\end{equation}
We point out some of basic properties of ${\cal F}$:

\begin{tabular}{l}
(a) The surface (\ref{equation:int-f}) has a simple singularity of type $E_7$ in the case $t_1=\cdots =t_5=t_7=0$.\\
(b) If the weights attached to  $x,y,z,t_j\>(j=1,2,3,4,5,7)$ are $2,
3,9/2,j\>(j=1,2,3,4,5,7)$, \\
respectively, then
$f_{\tau}$ is weighted homogeneous of total weight 9.\\
(c) If $t_7=0$, the origin of ${\bf C}^3$  is an isolated singular point of the surface  (\ref{equation:int-f}).\\
\end{tabular}

Since it is known that the degrees of basic invariants of ST34 are $6,12,18,24,30,42$,
 the set of weights of $t_1,\cdots,t_5,t_7$ and that of degrees of basic invariants are proportional.
This observation leads us to expect that $t_j$'s are regarded as  a certain set of basic invariants of ST34
and
 is actually true in the following sense.
 There is a polynomial $\delta(\tau)$ of $\tau$ such that  the surface defined by $f_{\tau}=0$ is smooth
 if and only if $t_7\cdot \delta\not=0$.
 If $\delta_{{\rm ST34}}$ is the discriminant of ST34,
 then it is possible to take $\delta_{\rm ST34}$ as the polynomial $\delta$.
 This leads to a simple and reasonable description of $\delta_{\rm ST34}$ in terms of basic invariants.
This is also regarded as the first step of an answer to the Arnold's problem for the group ST34.

So far we discussed the Arnol'd problem for the group ST34.
In this paper, we also introduce a  family $\tilde{\cal F}$ which 
 consists of the surfaces
${\tilde f}_{(s_3,\tau)}(x,y,z)=0$,
where
\begin{equation}
\label{equation:int-tilde-f}
{\tilde f}_{(s_3,\tau)}(x,y,z)=
y^3+x^3y+t_7x+t_5x^2+t_3x^3+t_1x^4+y(t_4x+t_2x^2)+s_3y^2-z^2.
\end{equation}
The members of ${\tilde{\cal F}}$ are also deformations of the simple singularity of type $E_7$.
It is noted here that $f_{\tau}$ is a special case of ${\tilde f}_{(s_3,\tau)}$.
The author constructed some families of (non-versal) deformations of simple singularities
related with the reflection groups of types $E_6,E_7,E_8$ and $H_4$ in \cite{Se3}.
The family $\tilde{\cal F}$  is one of them introduced there
in connection with an  algebraic potential of Frobenius manifold of type $E_7$.
We study the family ${\tilde{\cal F}}$  in the main text.

We now explain the the purposes of this paper.
Concerning the family $\tilde{\cal F}$,
we shall construct a polynomial ${\tilde\delta}(s_3,\tau)$ such that a member of
 $\tilde{\cal F}$ is singular  if and only if $t_7\cdot {\tilde\delta}(s_3,\tau)=0$
 and compare  the adjacent relation of   $\tilde{\cal F}$ with that of
 the versal deformation of the simple singularity of type $E_7$.
Concerning the family ${\cal F}$,
we determine 1-parameter subfamilies of the family  ${\cal F}$
corresponding to
 corank one reflection subgroups of  ST34.
To explain the result of this paper, we need some preparations.
We take the generators $P_1,P_2,P_3,Q_1,R_1,R_2$ of ST34 given in \cite{ST}, p.298.
For simplicity, we put
\begin{equation}
\label{equation:generator-ST34}
\rho_1=P_1,\rho_2=P_2,\rho_3=P_3,\rho_4=Q_1,\rho_5=R_1,\rho_6=R_2.
\end{equation}
Let  ${\rm ST34}_{j}$ denote the group generated by five generators $\rho_k\>(1\le k\le 6,\>k\not=j)$.
Then there is a curve $C_j$ of the parameter  space ${\bf C}^6$ (whose coordinate is $(t_1,\cdots,t_5,t_7)$)
passing through the origin
such that 
the surfaces corresponding to the parameters in $C_j-\{0,0,0,0,0,0\}$ are equisingular.
The purpose of this paper is to determine the type of simple singularities for the surfaces in question.
We denote its type by $X_{{\rm ST34}_j}$. 
Then  $X_{{\rm ST34}_j}$ is given
in the following table.
$$
\begin{array}{l|l|l|l}\hline
j&{\rm reflection}&{\rm ST34}_j&X_{{\rm ST34}_j}\\\hline\hline
1&\rho_1&W(A_3)\times W(A_2)&A_3+A_2\\\hline
2&\rho_2&G(3,3,4)\times W(A_1)&D_5+A_1\\\hline
3&\rho_3&ST33&E_6\\\hline
4&\rho_4&W(A_5)&A_5\\\hline
5&\rho_5&W(A_4)\times W(A_1)&A_4+A_1\\\hline
6&\rho_6&G(3,3,5)&D_6\\\hline
\end{array}
$$
(The group ST33 is the complex reflection group numbered as 33 in the list of \cite{ST}.)
We need an explanation to define the type $X_H$ for a complex reflection group $H$.
$$
\begin{array}{l|l}\hline
H&X_H\\\hline\hline
W(A_k)&A_k\\\hline
G(3,3,k)&D_{k+1}\\\hline
ST33&E_6\\\hline
\end{array}
$$
Moreover if $H=H_1\times H_2$ for some complex reflection groups $H_1,H_2$,
we write $X_H=X_{H_1}+X_{H_2}$.

This paper is organized as follows.
Section 2 starts with reviewing the root system of type $E_7$, construction of
 the generators of the ring of 
invariant polynomials, the definition of the  family of
semi-universal deformations of a simple singularity of type $E_7$.
At the last, we  explain the adjacent relation of the family in question.
In section 3, we introduce a family $\tilde{\cal F}$  of surfaces of ${\bf C}^3$ which is a deformation
of a simple singularity of type $E_7$ but not 
the same as the versal family.
The family $\tilde{\cal F}$ contains seven parameters $(s_3,t_1,\cdots,t_5,t_7)$.
We compute the discriminant of $\tilde{\cal F}$ and
confirm that each surface with corank one simple singularity appears in the family $\tilde{\cal F}$.
In section 4, we consider the family ${\cal F}$ which is obtained from $\tilde{\cal F}$
by substituting $s_3=0$.
Accordingly $\cal F$ contains six parameters   $(t_1,\cdots,t_5,t_7)$.
One of the reasons why we focus our attention on $\cal F$ is that the parameters
 $t_1,\cdots,t_5,t_7$ are regarded as basic invariant polynomials by the action of ST34.
 The precise meaning of this identification is explained in \S4.2.
In \S4.3, we study the adjacency relation between the surfaces contained in $\cal F$ with simple singularities
and the surface with an isolated singular point which is a simple singularity of type $E_7$.
Section 5 is devoted first to formulate an algebraic potentials,
next to construct the algebraic potential related with the family $\tilde{\cal F}$.

Finally we mention that the software Mathematica is used to obtain the results of this paper.


\section{A simple singularity of type $E_7$}
\label{section:two}
In this section, we start with reviewing the results on the root system of type $E_7$, the simple singularity of type $E_7$.
Then we compute 
 1-parameter deformations of simple singularity of type $E_7$.

\subsection{The root system of type $E_7$}
\label{subsection:2-1}
We first explain the definition and elementary properties of root systems
which will be used in the subsequent arguments, following \cite{Bour}.

Let ${\bf R}^8$ be an 8-dimensional Euclidean space with a basis $\{e_1,e_2,\cdots,e_8\}$,
an inner product $\langle\cdot,\cdot\rangle$ such that $\langle e_j,e_k\rangle=\delta_{jk}$.
Let $V$ be the hyperplane of ${\bf R}^8$ orthogonal to the vector $e_7+e_8$.
The root system of type $E_7$ denoted $R$ in this paper consists of 
vectors (called roots)  below:
$$
\begin{array}{l}
\pm e_i\pm e_j\>(1\le i<j\le 6),\quad \pm(e_7-e_8),\\
\pm\frac{1}{2}(e_7-e_8+\sum_{i=1}^6(-1)^{\nu(i)}e_i)
\>(\sum_{i=1}^6\nu(i):\>{\rm even})\\
\end{array}
$$
The roots
$$
\alpha_1=\frac{1}{2}(e_1+e_8)-\frac{1}{2}(e_3+e_4+e_5+e_6+e_7),
\alpha_2=e_1+e_2,\>\alpha_j=e_{j-1}-e_{j-2}\>(j=3,4,\cdots,7)
$$
form a system of simple roots.
Let $W(E_7)$ denote the group generated by the reflections with respect to the roots of $R$.

We now mention the invariant polynomials under $W(E_7)$  on $V_c=V\otimes_{\bf R}{\bf C}$.
Any vector of ${\tt v}\in V_c$ is
written as 
\begin{equation}
\label{equation:vector-v}
{\tt v}=\sum_{j=1}^6\xi_je_j+\xi_7(e_7-e_8).
\end{equation}
Let  $S={\bf C}[\xi_1,\xi_2,\cdots,\xi_7]$ be the polynomial ring over $V_c$.
It is known that the ring $S^{W(E_7)}$ of the $W(E_7)$-invariants of $S$ is
also a polynomial ring of 7 generators  which are homogeneous
and their homogeneous degrees are $2,6,8,10,12,14,18$.

We consider the weights of the 56-dimensional representation of the simple Lie algebra of type $E_7$
for later use.
Its highest weight  is
$\varpi_7=e_6+\frac{1}{2}(e_8-e_7)$.
Let $W(\varpi_7)=W(E_7)\cdot \varpi_7$ be the set of weights of the  56-dimensional representation.
We now represent the elements of $W(\varpi_7)$  explicitly.
We first define
$$
\Lambda_j=e_j+\frac{1}{2}(e_8-e_7)\>(1\le j\le 6),\quad \Lambda_7=e_1+\cdots+e_6
$$
and using these vectors, define
$$
\Lambda_{ij}(=\Lambda_{ji})=\Lambda_i+\Lambda_j-\frac{1}{3}\sum_{k=1}^7\Lambda_k\>(1\le i<j\le 7).
$$
Then $W(\varpi_7)$ consists of $\pm\Lambda_j,\>\pm\Lambda_{jk}$.
For the vector ${\tt v}\in V_c$ of the form (\ref{equation:vector-v}), we put
$$
u_j=\langle \Lambda_j,{\tt v}\rangle\>(1\le j\le 7),
\quad
u_{jk}=\langle \Lambda_{jk},{\tt v}\rangle\>(1\le j<k\le 7).
$$
We write $u_j\>(8\le j\le 28)$ instead of  $u_{jk}\>(1\le j<k\le 7)$
in an appropriate order.
It follows from the definition that $(u_1,\cdots,u_7)$ is regarded as a linear coordinate of $V_c$
and that $u_j\>(8\le j\le 28)$ are  expressed as linear combinations of $u_1,\cdots,u_7$.
We introduce
the polynomial
$$
\Psi_{{\tt v}}(X)=\prod_{i=1}^{28}(X^2-u_i^2).
$$
Let $\varepsilon_{\nu}$
be the $\nu$-th elementary symmetric polynomial of $u_1^2,\cdots, u_{28}^2$.
Obviously we have
$$
\Psi_{\tt v}(X)=X^{56}+\sum_{\nu=1}^{28}(-1)^{\nu}\varepsilon_{\nu}X^{56-2\nu}
$$
and the coefficients $\varepsilon_{\nu}$ are invariant under $W(E_7)$ as polynomials of $u_1,\cdots,u_7$.

\subsection{Semi-universal deformations of a simple singularity of type $E_7$}
\label{subsection:2-2}

We start with reviewing briefly the arguments  and results of \S9 in Shioda \cite{Sh}
which are useful 
to construct basic $W(E_7)$-invariants.

First of all we introduce the polynomial 
$$
f_{E_7}(x,y,z)=-z^2+y^3+y(p_0+p_1x+x^3)+q_0+q_1x+q_2x^2+q_3x^3+q_4x^4
$$
with parameters $p_0,p_1,q_0,q_1,q_2,q_3,q_4$.
The family of hypersurfaces
$f_{E_7}=0$ of ${\bf C}^3$ with the coordinate $(x,y,z)$ defines  semi-universal deformations
of a simple singularity of type $E_7$ which is a hypersurface of the form
$$
x^3y+y^3-z^2=0.
$$
For simplicity, we put
$\lambda=(p_0,p_1,q_0,q_1,q_2,q_3,q_4)$.
The computation of the construction of invariants of $W(E_7)$  in \cite{Sh} starts with the
determination of
curves of the form
\begin{equation}
\label{equation:lineareq}
y=ax+b,\>z=cx^2+dx+e\quad
(a,b,\cdots,e\in {\bf C})
\end{equation}
on the surface $f_{E_7}=0$.
The elimination of $y,z$ by (\ref{equation:lineareq}) implies
 the equation
\begin{equation}
\label{equation:fe7x}
f_{E_7}(x,ax+b,cx^2+dx+e)=0
\end{equation}
for $x$.
Then  condition that  (\ref{equation:fe7x})  holds for any $x$ is equivalent to
the system of equations for  $a,b,\cdots,e\in {\bf C}$ defined by
\begin{equation}
\label{equation:abcde}
\left\{
\begin{array}
{rcl}
c^2&=&a+q_4,\\
2cd&=&a^3+b+q_3,\\
d^2+2ce&=&3a^2b+p_1a+q_2,\\
2de&=&3ab^2+p_0a+p_1b+q_1,\\
e^2&=&b^3+p_0b+q_0.\\
\end{array}
\right.
\end{equation}
In virtue of  (\ref{equation:abcde}), it is easy to show that $a,b,e$ are expressed as rational
functions of $c,d$ over ${\bf Z}[p_0,p_1,q_0,q_1,q_2,q_3,q_4]$.
Moreover it follows that
$d\in {\bf Z}[p_0,p_1,q_0,q_1,q_2,q_3,q_4][c]$
and
$$
R_1(d)=d^3+\cdots=0,\quad R_2(d)=d^4+\cdots=0.
$$
Taking the resultant of  $R_1(d),\>R_2(d)$ with respect to $d$,
we obtain the algebraic equation for $c$
of the form
$$
{\tilde \Phi}(c)=c^{56}-36q_4c^{54}+594q_4^2c^{52}+(72q_3-6084q_4^3)c^{50}+\cdots=0,
$$
where the coefficients are contained in ${\bf Z}[p_0,p_1,q_0,q_1,q_2,q_3,q_4]$.
If $\lambda$ is generic, the equation ${\tilde \Phi}(c)=0$ has 56 solutions which are denoted by
$c_i\>(1\le i\le 56)$.
We note here  that ${\tilde \Phi}(-c)={\tilde \Phi}(c)$.
Since ${\tilde \Phi}$ depends on $\lambda$, we write $\Phi(X,\lambda)$ instead of ${\tilde \Phi}$:
$$
\Phi(X,\lambda)=\prod_{j=1}^{56}(X-c_j).
$$

One of the important statements in \cite{Sh} is the coincidence of $\Phi(X,\lambda)$ and $\Psi_{\tt v}(X)$,
namely the identity equation 
\begin{equation}
\label{equation:phipsi}
\Phi(X,\lambda)=\Psi_{\tt v}(X)
\end{equation}
holds if $p_0,p_1,q_0,q_1,q_2,q_4$ are polynomials of ${\tt v}$ to be determined below.
The comparison of the coefficients of 
the both sides of (\ref{equation:phipsi}) shows that
\begin{equation}
\label{equation:pqepsilon}
\left\{
\begin{array}{lll}
p_0 &=& \frac{1}{304749527040} (-318565 \varepsilon_1^6 + 6959520 \varepsilon_1^3 \varepsilon_3 + 
      58786560 \varepsilon_3^2 - 99283968 \varepsilon_1^2 \varepsilon_4 \\
      &&+ 
      321739776 \varepsilon_1 \varepsilon_5 
      - 
      564350976 \varepsilon_6), \\
     p_1 &=& \frac{1}{11197440} (12025 \varepsilon_1^4 - 129600 \varepsilon_1 \varepsilon_3 + 
      186624 \varepsilon_4), \\
     q_0 &=& (39246404257775 \varepsilon_1^9 - 876440477901600 \varepsilon_1^6 \varepsilon_3 + 
      4631876791257600 \varepsilon_1^3 \varepsilon_3^2 \\
      &&+ 1852115250585600 \varepsilon_3^3 + 
      1750210811020800 \varepsilon_1^5 \varepsilon_4 - 
            18355359088558080 \varepsilon_1^2 \varepsilon_3 \varepsilon_4\\
            && + 
      11416404882161664 \varepsilon_1 \varepsilon_4^2 - 1692182935434240 \varepsilon_1^4 \varepsilon_5 + 
      18973439502336000 \varepsilon_1 \varepsilon_3 \varepsilon_5\\
      && - 
            10826792911503360 \varepsilon_4 \varepsilon_5 
            - 
      336175411138560 \varepsilon_1^3 \varepsilon_6 - 14682832212787200 \varepsilon_3 \varepsilon_6 \\
      &&+ 
      6427167525273600 \varepsilon_1^2 \varepsilon_7 - 65167638862233600 \varepsilon_9)/
         1922184675880442265600, \\
  q_1 &=& \frac{1}{6999487137054720} (-4309074055 \varepsilon_1^7 + 
      79908139440 \varepsilon_1^4 \varepsilon_3 - 356716846080 \varepsilon_1 \varepsilon_3^2\\
      && - 
      73460991744 \varepsilon_1^3 \varepsilon_4 + 
            462203449344 \varepsilon_3 \varepsilon_4 + 143740790784 \varepsilon_1^2 \varepsilon_5 
            - 
      747200692224 \varepsilon_1 \varepsilon_6 \\
      &&+ 1828497162240 \varepsilon_7), \\
     q_2 &=& \frac{1}{470292480} (3575 \varepsilon_1^5 - 86400 \varepsilon_1^2 \varepsilon_3 + 
      435456 \varepsilon_1 \varepsilon_4 - 933120 \varepsilon_5), \\
  q_3 &=& (-169 \varepsilon_1^3 + 1296 \varepsilon_3)/93312,\\
   q_4 &=& -(\varepsilon_1/36).\\
\end{array}
\right.
\end{equation}
(The above equalities were given in \cite{Sh2}, Th.11.)
Hence we have
$$
{\bf Q}[p_0,p_1,q_0,q_1,q_2,q_3,q_4]=
{\bf Q}[\varepsilon_1,\varepsilon_3,\varepsilon_4,\varepsilon_{5},\varepsilon_{6},\varepsilon_{7},\varepsilon_{9}]
$$
which combined with the fact (cf. \S\ref{subsection:2-1})
 that $2,6,8,10,12,14,18$ are the degrees of basic $W(E_7)$-invariants of  $S^{W(E_7)}$
shows that  both
$\varepsilon_1,\varepsilon_3,\varepsilon_4,\varepsilon_{5},\varepsilon_{6},\varepsilon_{7},\varepsilon_{9}$
and $p_0,p_1,q_0,q_1,q_2,q_3,q_4$
form the basic invariants of  $S^{W(E_7)}$.

\subsection{Adjacency relation for the family $f_{E_7}=0$}
\label{subsection:2-3}
Adjacency relation of the simple singularity of type $E_7$ is already established
(for example, see Lamotke \cite{La} and Slodowy \cite{Sd1}, \S6.5).
The purpose of this subsection is to describe the details of the adjacency relation
from the view point of root systems.
Take a vector ${\tt v}\in V_c$ of the form 
(\ref{equation:vector-v})
and let $\alpha_i\>(1\le i\le 7)$ be simple roots of $\Delta(E_7)$ defined before.
Let $\varpi_j$ be the $j$-th weight, namely, $\varpi_j$
is the vector of $V$ defined by
$\alpha_k(\varpi_j)=\delta_{jk}$,
where $\delta_{jk}$ is the Kronecker's delta.
Using $\varpi_j$, we put, for $\xi\in{\bf C}$,

$$
\Psi_{\xi\varpi_i}(X)=\prod_{j=1}^7(X^2-\langle\Lambda_j,\xi\varpi_i\rangle^2)\cdot
\prod_{j=1}^6(\prod_{k=j+1}^7(X^2-\langle\Lambda_{jk},\xi\varpi_i\rangle^2))
$$
Then by the equation $\Phi(X,\lambda)=\Psi_{\xi\varpi_i}(X)$, 
it is possible to compute the value of each $\varepsilon_{\nu}$ and 
$p_0,p_1,q_0,q_1,q_2,q_3,q_4$ at ${\tt v}=\xi\varpi_i$.
We now put
$$
\lambda^{(i)}(\xi)=\lambda|_{{\tt v}=\xi\varpi_i}.
$$
Noting that the homogeneity of 
$p_0,p_1,q_0,q_1,q_2,q_3,q_4$, we have
$$
\lambda^{(i)}(\xi)=(\xi^{12}p_0^{(i)},\xi^8p_1^{(i)},\xi^{18}q_0^{(i)},\xi^{14}q_1^{(i)},\xi^{10}q_2^{(i)},\xi^6q_3^{(i)},\xi^2 q_4^{(i)}),
$$
where $(p_0^{(i)},p_1^{(i)},q_0^{(i)},q_1^{(i)},q_2^{(i)},q_3^{(i)}, q_4^{(i)})=\lambda|_{{\tt v}=\varpi_i}$.

The purpose of this subsection is to determine $\lambda^{(i)}(\xi)$
and the 1-parameter deformation of a simple singularity of type $E_7$  corresponding to $\lambda^{(i)}(\xi)$ 
for $i=1,2,\cdots,7$.
This is accomplished by (\ref{equation:pqepsilon}) and the result is given as follows.

\vspace{5mm}
{\bf Case $E_7$(1)}:$\alpha_1$

In this case, $2\xi\varpi_1=(0,0,0,0,0,0,\xi,-\xi)$ and 
$\Psi_{2\xi\varpi_1}(X)=X^{32}(X-2\xi)^{12}(X+2\xi)^{12}$.
Putting $\eta=\frac{4}{3}\xi^2$, we obtain
\begin{equation}
\label{equation:vece7-1}
\lambda^{(1)}=( -2 \eta^6,  -3 \eta^4,  -\eta^9,  -2 \eta^7, 0, 
 2 \eta^3,  \eta)
\end{equation}
 and
$$-\eta^9 - 2 \eta^7 x + 2 \eta^3 x^3 + \eta x^4 - 2 \eta^6 y - 3 \eta^4 x y + 
   x^3 y + y^3 - z^2=0.$$
By the change of variables $(x,y,z)=(x_1-\eta^2,y_1-\eta x_1,z)$, this equation turns out to be
\begin{equation}
\label{equation:alpha-1case}
f_{\alpha_1}=x_1^3y_1+y_1^3-3\eta x_1y_1^2-z^2=0.
\end{equation}
In virtue of \cite{La}, p.212. we find that if $\eta\not=0$, this surface in the  $(x_1,y_1,z)$-space has a simple singularity of type $D_6$ at the origin.

\vspace{5mm}
{\bf Case $E_7$(2)}: $\alpha_2$

In this case, $2\xi\varpi_2=(\xi,\xi,\xi,\xi,\xi,\xi,-2\xi,2\xi)$
and
$\Psi_{2\xi\varpi_2}(X)=(X^2-36\xi^2)^7(X^2-4\xi^2)^{21}$.
Putting $\eta=8\xi^2$, 
and by   the change of variables $(x,y,z)=(x_1+\frac{1}{4}\eta^2,y_1+\frac{7}{6}\eta x_1+\frac{23}{24}\eta^3,z)$,
we obtain
\begin{equation}
\label{equation:vece7-2}
\lambda^{(2)}=( -812 \eta^6/729,  -28 \eta^4/27, 
 13384 \eta^9/19683, 4232 \eta^7/2187, 28 \eta^5/9, 
   70 \eta^3/27, 7 \eta/6)
 \end{equation}
   and
\begin{equation}
\label{equation:alpha-2case}
f_{\alpha_2}=\frac{8}{9} \eta^5 x_1^2- \frac{52}{27}  +\frac{7}{6} \eta x_1^4 + 
  \frac{8}{3} \eta^4 x_1 y_1 -\frac{10}{3} \eta^2 x_1^2 y_1 + x_1^3 y_1 + 
   2 \eta^3 y_1^2 + y_1^3 - z^2=0.
 \end{equation}  
If $\eta\not=0$, this surface in the  $(x_1,y_1,z)$-space has a simple singularity of type $A_6$ at the origin.
This is shown as follows.
Changing the variable $y_1$ with $y_2$ by
$$
y_1=y_2-\frac{2}{3}\eta x_1+\frac{x_1^2}{2\eta}+\frac{x_1^3}{4\eta^3}+\frac{x_1^4}{16\eta^5}-\frac{x_1^5}{8\eta^7}-\frac{7x_1^6}{32\eta^9},
$$
we find that
$$
f_{\alpha_2}=y_2^2\varphi(x_1,y_2)+x_1^7\psi(x_1,y_2)-z^2,
$$
where $\varphi$ and $\psi$ are polynomials of $x_1,y_2$ and that $\varphi(0,0)\not=0,\>\psi(0,0)\not=0$.
As a consequence,
$f_{\alpha_2}=0$ has a simple singularity of type $A_6$ at $(x_1,y_2,z)=(0,0,0)$.

\vspace{5mm}
   {\bf Case $E_7$(3)}: $\alpha_3$

   In this case, $3\xi\varpi_3=(-\xi,\xi,\xi,\xi,\xi,\xi,-3\xi,3\xi)$
   and
$\Psi_{3\xi\varpi_3}(X)=X^{20}(X-64\xi^2)^{6}(X^2-16\xi^2)^{12}$.
Putting $\eta=16\xi^2$, 
and by
the change of variables $(x,y,z)=(x,y_1+\frac{1}{3}\eta^3,z)$,
we obtain
\begin{equation}
\label{equation:vece7-3}
\lambda^{(3)}=( -\eta^6/3,  -\eta^4,  2 \eta^9/27,  \eta^7/3, 
 \eta^5,  5 \eta^3/3,  \eta)
\end{equation}
and
\begin{equation}
\label{equation:alpha-3case}
f_{\alpha_3}=\eta^5x^2+2\eta^3x^3+\eta x^4-\eta^4 xy_1+x^3 y_1+\eta^3 y_1^2+y_1^3-z^2=0.
\end{equation}
If $\eta\not=0$, this surface in the  $(x,y_1,z)$-space has a simple singularity of type $A_1$ at the origin
and a simple singularity of type $A_5$ at $(x,y_1,z)=(-\eta^2,0,0)$.

\vspace{5mm}
{\bf Case $E_7$(4)}: $\alpha_4$

In this case, $\xi\varpi_4=(0,0,\xi,\xi,\xi,\xi,-2\xi,2\xi)$
and
$\Psi_{\xi\varpi_4}(X)=X^{12}(X^2-36\xi^2)^4(X^2-16\xi^2)^{6}(X^2-4\xi^2)^{12}$.
Putting $\eta=4\xi^2$,  
we obtain
\begin{equation}
\label{equation:vece7-4}
\lambda^{(4)}=( -46 \eta^6/3, -11 \eta^4,  812 \eta^9/27, 
 146 \eta^7/3,  34 \eta^5,  38 \eta^3/3,  2 \eta)
\end{equation}
 and
\begin{equation}
\label{equation:alpha-4case}
f_{\alpha_4}=\frac{812}{27}\eta^9+\frac{146}{3}\eta^7x+34\eta^5x^2+\frac{38}{3}\eta^3x^3+2\eta x^4
-\frac{46}{3}\eta^6y-11\eta^4xy+x^3y+y^3-z^2=0.
\end{equation}
If $\eta\not=0$, this surface in the  $(x,y,z)$-space has  simple singularities  at three points
$(-\eta^2,\frac{4}{3}\eta^3,0),$ $(-2\eta^2,-\frac{2}{3}\eta^3,0),$
$(-5\eta^2,\frac{16}{3}\eta^3,0)$.
Their types are $A_2,\>A_1,\>A_3$, respectively.

\vspace{5mm}
{\bf Case $E_7$(5)}: $\alpha_5$

In this case, $\xi\varpi_5=(0,0,0,\xi,\xi,\xi,-\frac{3}{2}\xi,\frac{3}{2}\xi)$
and
$\Psi_{\xi\varpi_5}=(X^2-25\xi^2)^3(X^2-9\xi^2)^{10}(X^2-\xi^2)^{15}$.
Putting $\eta=2\xi^2$,  
we obtain
\begin{equation}
\label{equation:vece7-5}
\lambda^{(5)}=(  -100 \eta^6/3,   -20 \eta^4,   3400 \eta^9/27, 
  520 \eta^7/3,   92 \eta^5,  70 \eta^3/3, 
    5 \eta/2)
 \end{equation}
    and
\begin{equation}
\label{equation:alpha-5case}
f_{\alpha_5}=\frac{3400}{27}\eta^9+\frac{520}{3}\eta^7x+92\eta^5x^2+\frac{70}{3}\eta^3x^3+\frac{5}{2}\eta x^4+\frac{100}{3}\eta^6 x
-20\eta^4 xy+x^3y+y^3-z^2=0.
\end{equation}
If $\eta\not=0$, this surface in the  $(x,y,z)$-space has a  simple singularity of type $A_4$ at 
$(-2\eta^2,\frac{2}{3}\eta^3,0)$
and a simple singularity of type $A_2$ at $(-10\eta^2,\frac{50}{3}\eta^3,0)$.

\vspace{5mm}
{\bf Case $E_7$(6)}: $\alpha_6$

In this case, $\xi\varpi_6=(0,0,0,0,\xi,\xi,-\xi,\xi)$
and
$\Psi_{\xi\varpi_6}=X^{20}(X^2-16\xi^2)^2(X^2-4\xi^2)^{16}$.
Putting $\eta=\frac{4}{3}\xi^2$,  
by the change of variables $(x,y,z)=(x_1-\eta^2,y,z)$,
we obtain
\begin{equation}
\label{equation:vece7-6}
\lambda^{(6)}=( -2 \eta^6,  -3 \eta^4,  4 \eta^9,   14 \eta^7, 
   18 \eta^5,   10 \eta^3,  2 \eta)
\end{equation}
   and
\begin{equation}
\label{equation:alpha-6case}
f_{\alpha_6}=2\eta^3x_1^3+2\eta x_1^4-3\eta^2 x_1^2y+x_1^3y+y^3-z^2=0.
\end{equation}
If $\eta\not=0$, this surface in the  $(x_1,y,z)$-space has a  simple singularity of type $D_5$ at the origin
and a simple singularity of type $A_1$ at $(-9\eta^2,18\eta^3,0)$.
We only show that the surface  $f_{\alpha_6}=0$ has a simple singularity of type $D_5$ at $(x_1,y,z)=(0,0,0)$.
By the change of variables $x_1,y_1$ with $x_2,y_2$:
$$
\begin{array}{lll}
x_1&=&\frac{x_2}{3\eta}-\frac{y_2}{2\eta}-\frac{x_2^2}{54\eta^4},\\
y&=&\frac{x_2}{3}+y_2+\frac{x_2^2}{27\eta^3},\\
\end{array}
$$
we find that
$$
f_{\alpha_6}=x_2(\varphi(x_2,y_2)x_2^3+\psi(x_2,y_2)y_2^2)-z^2,
$$
where $\varphi,\>\psi$ are polynomials of $x_2,y_2$ such that $\varphi(0,0)\not=0,\>\psi(0,0)\not=0$.
Then it is easy to see that $f_{\alpha_6}=0$ has a simple singularity of type $D_5$ at $(x_1,y,z)=(0,0,0)$.

\vspace{5mm}
{\bf Case $E_7$(7)}: $\alpha_7$

In this case, $\xi\varpi_7=(0,0,0,0,0,\xi,-\frac{1}{2}\xi,\frac{1}{2}\xi)$
and
$\Psi_{\xi\varpi_7}=(X^2-9\xi^2)(X^2-\xi^2)^{27}$.
Putting $\eta=\xi^2$,  
we obtain
\begin{equation}
\label{equation:vece7-7}
\lambda^{(7)}=(0,0,0,0,0,0,\eta)
\end{equation}
and
\begin{equation}
\label{equation:alpha-7case}
f_{\alpha_7}=\eta x^4+x^3y+y^3-z^2=0.
\end{equation}
If $\eta\not=0$, this surface in the  $(x,y,z)$-space has a  simple singularity of type $E_6$ at the origin.
This follows from \cite{La}, p.212.

\begin{remark}
We regard  $f_{\alpha_j}(x,y,z)\>(1\le j\le 7)$ as polynomials of $(x,y,z)$ with a parameter $\eta$
in the above consideration.
Regarding $f_{\alpha_j}(x,y,0)$ as a polynomial of  $(\eta,x,y)$
and denoting it by $f_j(\eta,x,y)$, we observe that $f_j(\eta,x,y)=0$ defines a free divisor 
in ${\bf C}^3$ in the sense of K. Saito.
We now recall the definition of polynomials $F_{B,j}(x,y,z)\>(1\le j\le 7)$ in \cite{Se1}.
Then by changing variables, we find 
the existence of a 1-1 correspondence between $\{f_j(\eta,x,y)|\>1\le j\le 7\}$
and $\{F_{B,j}(x,y,z)|\>1\le j\le 7\}$ in the following manner:
$$
\begin{array}{l|l|l|l|l|l|l}\hline
f_1&f_2&f_3&f_4&f_5&f_6&f_7\\\hline
F_{B,3}&F_{B,7}&F_{B,2}&F_{B,1}&F_{B,6}&F_{B,4}&F_{B,5}\\\hline
\end{array}
$$

\end{remark}


\section{A new family of surfaces of ${\bf C}^3$ which is
a deformation of a simple singularity of type $E_7$}
\label{section:three}

In this section, we start with introducing a family of surfaces of ${\bf C}^3$ denoted by $\tilde{\cal F}$
and then compute the discriminant of $\tilde{\cal F}$.
Finally we confirm the appearance in $\tilde{\cal F}$  of  the surfaces corresponding to
corank one subdiagrams of the root system $R$.

\subsection{The new family of surfaces of ${\bf C}^3$}
\label{subsection:3-1}

We first introduce a polynomial
$
{\tilde f}_{(s_3,\tau)}(x,y,z)$ ($\tau=(t_1,t_2,t_3,t_4,t_5,t_7)$) defined by
\begin{equation}
\label{equation:polynomial-algebraicFrobenius}
{\tilde f}_{(s_3,\tau)}(x,y,z)=
y^3+x^3y+t_7x+t_5x^2+t_3x^3+t_1x^4+y(t_4x+t_2x^2)+s_3 y^2-z^2
\end{equation}
If the weights of $x,y,z,s_3,t_j\>(j=1,2,3,4,5,7)$ are $ 2,3,9/2,3,j$, respectively,
${\tilde f}_{(s_3\tau)}(x,y,z)$ is weighted homogeneous of total weight 9.

Using the polynomial ${\tilde f}_{(s_3,\tau)}$, we define
the surface ${\tilde S}(s_3,\tau)$ of the $(x,y,z)$-space by
$$
{\tilde S}(s_3,\tau):{\tilde f}_{(s_3,\tau)}(x,y,z)=0.
$$
Then
$$
{\tilde{\cal F}}=\{{\tilde S}(s_3,\tau)|\>(s_3,\tau)\in{\bf C}^7\}
$$
is a family of surfaces of ${\bf C}^3$
and ${\tilde S}(0,0,0,0,0,0,0)$ has a simple singularity of type $E_7$ at the origin of ${\bf C}^3$.

\begin{lemma}

(i)
If
\begin{equation}
\label{equation:pqst}
\left\{
\begin{array}{lll}
p_0 &=&\frac{1}{27} (-9 s_3^2 + 2 t_2^3 - 9 t_2 t_4),\\ 
  p_1 &=& \frac{1}{3} (-t_2^2 + 3 t_4),\\ 
     q_0 &=&\frac{1}{81} (6 s_3^3 - 2 s_3 t_2^3 + t_1 t_2^4 - 3 t_2^3 t_3 + 
      9 s_3 t_2 t_4 + 9 t_2^2 t_5 - 27 t_2 t_7),\\ 
     q_1 &=&\frac{1}{27} (3 s_3 t_2^2 - 4 t_1 t_2^3 + 9 t_2^2 t_3 - 9 s_3 t_4 - 
      18 t_2 t_5 + 27 t_7),\\ 
      q_2 &=&\frac{1}{3} (2 t_1 t_2^2 - 3 t_2 t_3 + 3 t_5),\\ 
     q_3 &=&\frac{1}{3} (-s_3 - 4 t_1 t_2 + 3 t_3),\\ 
     q_4 &=& t_1,
\\
\end{array}
\right.
\end{equation}
then
$$
{\tilde f}_{(s_3,\tau)}(x_1-\frac{1}{3}t_2,y_1-\frac{1}{3}s_3,z)
=
f_{E_7}(x_1,y_1,z).
$$

(ii) For any $(p_0,p_1,q_0,q_1,q_2,q_3,q_4)\in{\bf C}^7$,
there is $(s_3,\tau)\in{\bf C}^7$ satisfying the equation (\ref{equation:pqst}).

\end{lemma}

{\bf Proof:}
(i) follows from a direct computation.

(ii)
For given  $(p_0,p_1,q_0,q_1,q_2,q_3,q_4)\in{\bf C}^7$,
we now solve 
(\ref{equation:pqst}), being regarded as a system of algebraic equations for $(s_3,\tau)$.
It is clear that $t_1,t_4$ are given by
\begin{equation}
\label{equation:condition-t1t4}
\left\{
\begin{array}{lll}
t_1 &=& q_4, \\
  t_4 &=& \frac{1}{3} (3 p_1 + t_2^2).\\
     \end{array}
     \right.
      \end{equation}
      Then $s_3,t_2$ are solutions of the algebraic equations
      \begin{equation}
      \label{equation:s3t2}
      \left\{
      \begin{array}{l}
-81 q_0 + 6 s_3^3 - 27 q_1 t_2 - 9 q_2 t_2^2 - 3 q_3 t_2^3 - 
   q_4 t_2^4=0, \\
   27 p_0 + 9 s_3^2 + 9 p_1 t_2 + t_2^3=0.\\
   \end{array}
   \right.
   \end{equation}
Taking a solution  $(s_3,t_2)$ of   (\ref{equation:s3t2}), we have
$t_1,t_4$ by (\ref{equation:condition-t1t4}) and
$t_3,t_5,t_7$ by
\begin{equation}
\left\{
\begin{array}{lll}
t_3 &=& \frac{1}{3} (3 q_3 + s_3 + 4 q_4 t_2), \\
  t_5 &=&\frac{1}{3} (3 q_2 + 3 q_3 t_2 + s_3 t_2 + 2 q_4 t_2^2), \\
     t_7 &=& \frac{1}{27} (27 q_1 + 9 p_1 s_3 + 18 q_2 t_2 + 9 q_3 t_2^2 + 
      3 s_3 t_2^2 + 4 q_4 t_2^3).\\
      \end{array}
      \right.
      \end{equation}
In this way, all the solutions   $(s_3,\tau)$ to (\ref{equation:pqst}) are obtained.
We  note here that 
 (\ref{equation:s3t2}) implies that
\begin{equation}
s_3 = -\frac{3 (81 q_0 + 27 q_1 t_2 + 9 q_2 t_2^2 + 3 q_3 t_2^3 + 
          q_4 t_2^4)}{2 (27 p_0 + 9 p_1 t_2 + t_2^3)}.
          \end{equation}
        
 \subsection{The discriminant of the family ${\tilde{\cal F}}$
 }\label{subsection:3-2}

 Putting ${\tilde f}={\tilde f}_{(s_3,\tau)}$ for simplicity, we define
$$
{\tilde C}_{\tilde f}=\{(x,y,z;s_3,\tau)\in{\bf C}^{10}|\>{\tilde f}=\partial_x{\tilde f}=\partial_y{\tilde f}=\partial_z{\tilde f}=0\}
$$
and
$$
C_{\tilde f}=\{(s_3,\tau)\in{\bf C}^{7}|\>\exists (x,y,z)\in{\bf C}^3\>{\rm s.t.}\>(x,y,z;s_3,\tau)\in{\tilde C}_{\tilde f}\}.
$$
It is clear that if $t_7=0$, then ${\tilde f}(0,0,0)=\partial_x{\tilde f}(0,0,0)=\partial_y{\tilde f}(0,0,0)=0$.
This means that $C_{\tilde f}\subseteq \{(s_3,\tau)|\>t_7=0\}$.
It is provable that there is a polynomial  ${\tilde \delta}(s_3,\tau)$ such that
$$
C_{\tilde f}=\{(s_3,\tau)\in{\bf C}^7|\>t_7{\tilde \delta}(s_3,\tau)=0\}.
$$

We are going to construct $\tilde\delta$ concretely.
For this purpose, we  introduce polynomials of $x$
$$
\begin{array}{lll}
A_0&=&t_7x+t_5x^2+t_3x^3+t_1x^4,\\
A_1&=&t_4x+t_2x^2,\\
A_2&=&s_3,\\
M_0&=&a_1+a_2x+a_3x^2+a_4x^3+a_5x^4,\\
M_1&=&a_6+a_7x,\\
N_0&=&b_1x+b_2x^2,\\
N_1&=&b_3+b_4x+b_5x^2+b_6x^3,\\
N_2&=&b_7+b_8x+b_9x^2,\\
P_0&=&c_1+c_2x+c_3x^2+c_4x^3+c_5x^4+c_6x^5,\\
P_1&=&c_7+c_8x+c_9x^2+c_{10}x_3+c_{11}x^4,\\
P_2&=&c_{12}+c_{13}x,\\
\end{array}
$$
and
$$
\begin{array}{lll}
g_1&=&M_0+M_1y,\\
g_2&=&N_0+N_1y+N_2y^2,\\
g_3&=&P_0+P_1y+P_2y^2.\\
\end{array}
$$
Then $g_1,g_2,g_3$ are polynomials of $(x,y)$ and $a_i,b_j,c_k$ are constants to be determined.
Note that ${\tilde f}=A_0+A_1y+A_2y^2+y^3-z^2$.
Using $g_1,g_2,g_3$, we define
$$
H=g_1{\tilde f}+g_2\partial_x{\tilde f}+g_3\partial_y{\tilde f}-\frac{1}{2}g_1\partial_z{\tilde f}.
$$
Then $H$ is independent of $z$.
Putting ${\tilde g}={\tilde f}+z^2$, we find that
$$
H=g_1{\tilde g}+g_2\partial_x{\tilde g}+g_3\partial_y{\tilde g}.
$$

We expand $H$ as a polynomial of $y$, namely,
$$
H=H_0(x)+H_1(x)y+H_2(x)y^2+H_3(x)y^3+H_4(x)y^4.
$$
In particular,
$$
\begin{array}{lll}
H_2&=&3P_0+A_2M_0+A_1M_1+\cdots,\\
H_3&=&3P_1+M_0+A_2M_1\cdots,\\
H_4&=&3P_2+M_1.\\
\end{array}
$$
Accordingly  it is possible to choose $c_k$'s so that $H_2=H_3=H_4=0$ hold.
Needless to say, $c_k$'s are
expressed as linear combinations of  $a_i$'s and $b_j$'s with polynomial coefficients of $(s_3,\tau)$.
Then $H$ turns out to be
$$
H=K_0(x)+K_1(x)y,
$$
where $K_0(x),K_1(x)$ are polynomials of degree 8, 7, respectively
and their coefficients are linear combinations of $a_i$'s and $b_j$'s.
Taking into account of the coefficients of $x^m\>(5\le m\le 8)$ in $K_0$
 and those of  $x^m\>(3\le m\le 7)$ in $K_1$, we can eliminate $b_1,\cdots,b_9$
 under the condition that these coefficients are 0.
 As a result, $H$ turns out to be
 $$
 H=L_0(x)+L_1(x)y,
 $$
 where $L_0(x),L_1(x)$ are polynomials of degrees 4, 2, respectively
and their coefficients are linear combinations of $a_i$'s.
It is underlined here that all the coefficients of $L_0,L_1$ are polynomials of $(s_3,\tau)$,
that $H(0,0)=K_0(0)=0$
and that
\begin{equation}
\label{equation:mat-A}
L_0(x)+L_1(x)y|_{t_1=\cdots=t_5=s_3=0}=\frac{7}{9}t_7(a_1x+a_2x^2+a_3x^3+a_4x^4+y(a_5t_7+a_6x+a_7x^2)).
\end{equation}
These properties follow from the computation.
Then 
there are polynomials $A_{ij}(s_3,\tau)\>(1\le j\le 7)$ of $(s_3,\tau)$
such that
\begin{equation}
\label{equation:comp-H}
H=\sum_{i=1}^7a_i\left(\sum_{j=1}^4A_{ij}x^j+\sum_{k=1}^3A_{i,k+4}x^{k-1}y\right).
\end{equation}
We now put $A=(A_{ij})$ which is a $7\times 7$ matrix.
Then by a direct computation, we find that
\begin{equation}
\label{equation:matrix-B}
A|_{t_1=t_2=t_3=t_4=0}=
\frac{7}{9}\left(
\begin{array}{ccccccc}
t_7&\frac{5}{7}t_5&-\frac{1}{7}s_3&0&-\frac{2}{7}s_3^2&0&0\\
0&t_7&\frac{5}{7}t_5&-\frac{1}{7}s_3&0&-\frac{2}{7}s_3^2&0\\
0&0&t_7&\frac{5}{7}t_5&-\frac{1}{7}s_3t_7&-\frac{2}{7}s_3t_5&0\\
0&0&0&t_5&\frac{5}{7}t_5t_7&\frac{10t_5^2-s_3t_7}{7}&-\frac{12}{7}s_3t_5\\
\frac{4}{7}s_3t_5t_7&\frac{8}{7}s_3t_5^2&0&0&t_7^2&\frac{19}{7}t_5t_7&\frac{5(2t_5^2-3s_3t_7)}{7}\\
\frac{1}{21}s_3t_7&\frac{2}{21}s_3t_5&\frac{2}{21}s_3^2&0&\frac{4}{21}s_3^3&t_7&\frac{5}{7}t_5\\
-\frac{5}{21}t_5t_7&\frac{-10t_5^2+s_3t_7}{21}&\frac{2}{21}s_3t_5&\frac{2}{21}s_3^2&0&\frac{4}{21}s_3^3&t_7\\
\end{array}
\right).
\end{equation}
As a consequence,  it follows that
$A$ is generically invertible.

\begin{theorem}
\label{theorem:discriminant1}
There is a polynomial $\tilde{\delta}$ of $(s_3,\tau)$
with the following properties.

(i) There is a non-zero constant $k_0$ such that
$$
\det A=k_0t_7{\tilde\delta}.
$$

(ii) 
As a polynomial of $t_7$,  $\tilde{\delta}$ is expressed as
$$
{\tilde\delta}(s_3,\tau)=t_7^7+\sum_{j=1}^7{\tilde A}_j(s_3,t_1,\cdots,t_5)t_7^{7-j}.
$$

(iii)
${\tilde{\delta}}$ is irreducible as a polynomial of  $(s_3,\tau)$.

(iv) $C_{\tilde f}=\{(s_3,\tau)\in {\bf C}^7|\>t_7{\tilde\delta}=0\}$

(v) There is a polynomial $\delta_{\rm ST34}$ of $\tau$
such that
$$
{\tilde\delta}(0,\tau)=t_7\delta_{\rm ST34}(\tau).
$$

\end{theorem}

{\bf Proof.}
(i), (ii), (v) are shown by a direct computation.

(iii) To prove (iii), we put $\tilde{\delta}_0(s_3,t_5,t_7)=\tilde{\delta}(s_3,0,0,0,0,t_5,t_7)$.
If $\tilde{\delta}$ is reducible, so is $\tilde{\delta}_0$.
Therefore it suffices to prove that   $\tilde{\delta}_0$ is irreducible as a polynomial of $(s_3,t_5,t_7)$.
If $\tilde{\delta}_0$ is reducible,
 its factor is one of
$t_7,\>t_7^2+a_1m_1t_5s_3^3,\>t_7^3+n_1t_5^3s_3^2+n_2s_3^7$ for some constants $m_1,n_1,n_2$.
It is easy to check that all of $t_7,\>t_7^2+a_1m_1t_5s_3^3,\>t_7^3+n_1t_5^3s_3^2+n_2s_3^7$ 
are not factors of $\tilde{\delta}_0$.
Then $\tilde{\delta}_0$ is irreducible, so is $\tilde{\delta}$.

(iv) 
We put 
$$D_{\tilde f}^{(1)}=\{(s_3,\tau)\in{\bf C}^7|\>{\tilde\delta}(s_3,\tau)=0\},\quad
D_{\tilde f}^{(2)}=\{(s_3,\tau)\in{\bf C}^7|\>t_7=0\}
$$
and
$$D_{\tilde f}=D_{\tilde f}^{(1)}\cup D_{\tilde f}^{(2)}
$$
 for simplicity.

We first prove $
C_{\tilde f}\subseteq D_{\tilde f}$.
For this purpose, we take a point  $(s_3,\tau)\in C_{\tilde f}$.
Then there is $(u,v,w)\in {\bf C}^3$ such that   $(u,v,w;s_3,\tau)\in {\tilde C}_{\tilde f}$.
If $(u,v,w)=(0,0,0)$, then $t_7=0$.
Next we consider the case where $(u,v,w)\not=(0,0,0)$.
Since $\partial_z{\tilde f}=-2z$, it follows that $w=0$.
Noting this, we may assume $w=0$ from the first.
We consider ${\tt u}=(u,u^2,u^3,u^4,v,uv,u^2v)$ which is not a zero vector.
In virtue of  (\ref{equation:comp-H}), we find that
$$
\sum_{i=1}^7a_i\left(\sum_{j=1}^4A_{ij}u^j+\sum_{k=1}^3A_{i,k+4}u^{k-1}v\right)=0
$$
for any $a_1,\cdots,a_7$.
Then
$A{}^t{\tt u}=O$.
This means  that  $\det A=0$ for any $(s_3,\tau)\in C_{\tilde f}$
with the condition that there is $(u,v,0;s_3,\tau)\in {\tilde{\cal C}}_{\tilde f},\>(u,v)\not=(0,0)$.
As a consequence, we find that
$
C_{\tilde f}\subseteq  D_{\tilde f}
$.

We next prove
$D_{\tilde f}\subseteq C_{\tilde f}$.
Since it is clear that $D_{\tilde f}^{(2)}\subset C_{\tilde f}$,
it suffices to prove $D_{\tilde f}^{(1)}\subset C_{\tilde f}$,
For this purpose, we put
$$
D_{\tilde f,\>{\rm smooth}}^{(1)}=D_{\tilde f}-(D_{{\tilde f},t_1}\cap \cdots \cap D_{{\tilde f},t_5}\cap
D_{{\tilde f},t_7}\cap D_{{\tilde f},s_3}),
$$
where 
$$
D_{{\tilde f},t_j}=\{(s_3,\tau)\in D_{\tilde f}|\>\partial_{t_j}{\tilde\delta}=0\},
\quad
D_{{\tilde f},s_3}=\{(s_3,\tau)\in D_{\tilde f}|\>\partial_{s_3}{\tilde\delta}=0\}.
$$
Then $D_{\tilde f,\>{\rm smooth}}^{(1)}$ consists of smooth point of  $D_{\tilde f}^{(1)}$.
Since $\tilde\delta$ is irreducible,
it is sufficient to prove that for a point $(s_3,\tau)\in  D_{\tilde f,\>{\rm smooth}}^{(1)}$,
there is $(u,v)\in {\bf C}^2$ such that $(u,v,0;s_3,\tau)\in {\tilde C}_{\tilde f}$.
We put $B=A_{t_1=\cdots=t_4=0}$.
Then, noting
(\ref{equation:matrix-B}),
we find that
$$
\det(B)=\frac{7^7}{9^7}
\left(t_7^7-\frac{225}{343}s_3^3t_5t_7^5-\frac{25}{823543}(3375s_3^5-686t_5^3)s_3^2t_7^4+\cdots\right).
$$
Since $\det(B)$ coincides with  ${\tilde\delta}_0(s_3,t_5,t_7)$ introduced in (iii) up to a constant factor,
it is easy to compute $\partial_{t_j}{\tilde\delta}_0\>(j=5,7)$ and $\partial_{s_3}{\tilde\delta}_0$.
As a consequence, we conclude that
 $D_{{\tilde f},{\rm smooth}}\cap \{(s_3,\tau)\in{\bf C}^7|\>t_1=\cdots=t_4=0\}\not=\emptyset$.

We now take $(s_3^{\ast},\tau^{\ast})\in D_{{\tilde f},{\rm smooth}}$ of the form $\tau^{\ast}=(0,0,0,0,t_5^{\ast},t_7^{\ast})$
and fix it for a moment. 
If we prove  $(s_3^{\ast},\tau^{\ast})\in C_{{\tilde f}}$, we conclude 
$D_{\tilde f}\subseteq C_{\tilde f}$
 by the reason that  $(s_3^{\ast},\tau^{\ast})$ is a generic point of $D_{{\tilde f},{\rm smooth}}$.
 We are going to prove  $(s_3^{\ast},\tau^{\ast})\in C_{{\tilde f}}$.
By definition,  ${\tilde\delta}(s_3^{\ast},\tau^{\ast})={\tilde\delta}_0(s_3^{\ast},t_5^{\ast},t_7^{\ast})=0$.
Let $B^{\ast}$ be the matrix $B_{(s_3,t_5,t_7)=(s_3^{\ast},t_5^{\ast},t_7^{\ast})}$.
Since $\det(B^{\ast})=0$, there is a non-zero vector ${\tt u}=(u_1,u_2,\cdots,u_7)$ such that
$B^{\ast}\cdot{}^t {\tt u}=\vec{0}$.

We now treat $B$ instead of $B^{\ast}$
and define $h_j=h_j(s_3,t_5,t_7)$ by ${}^t(h_1,\cdots,h_7)=B\cdot{}^t{\tt u}$.
By solving the system of equations
$$
h_j=0\>(j=2,\cdots,7)
$$
for $u_2,\cdots,u_7$, we find the existence of  rational functions  $\gamma_j=\gamma_j(s_3,t_5,t_7)\>(j=2,\cdots,7)$ of $s_3,t_5,t_7$
such that $u_j=\gamma_j u_1\>(j=2,\cdots,7)$.
Then, the vector $\tilde{\tt u}=(1,\gamma_2,\cdots,\gamma_7)$ satisfies that $B\cdot{}^t{\tilde{\tt u}}=
{}^t(\gamma_0\tilde\delta_0,0,0,0,0,0,0)$ for a non-zero rational function$\gamma_0$ of $s_3,t_5,t_7$.
We now observe
$\gamma_2\gamma_j\equiv \gamma_{j+1}\>({\rm mod}\>\> \tilde\delta_0)\>(j=2,3,5,6)$.
As a consequence,
$$
(1,\gamma_2,\cdots,\gamma_7)\equiv (1,\gamma_2,\gamma_2^2,\gamma_2^3,\gamma_5,\gamma_2\gamma_5,\gamma_2^2\gamma_5)\>({\rm mod}\>\tilde\delta_0).
$$
Noting this, we put $u=\gamma_2,v=\gamma_2\gamma_5$.
Returning back to the case $(s_3,t_5,t_7)=(s_3^{\ast},t_5^{\ast},t_7^{\ast})$,
we put $\gamma_j^{\ast}=\gamma_j(s_3^{\ast},t_5^{\ast},t_7^{\ast})\>(j=2,5),
\>u^{\ast}=\gamma_2^{\ast},v^{\ast}=\gamma_2^{\ast}\gamma_5^{\ast}$.
Then the vector
${\tt u}^{\ast}=(u^{\ast},(u^{\ast})^2,(u^{\ast})^3,(u^{\ast})^4,v^{\ast},u^{\ast}v^{\ast},(u^{\ast})^2v^{\ast})$
satisfies that $B^{\ast}\cdot{}^t{\tt u}^{\ast}={}^t\vec{0}$.

By a direct computation, we may take as $u,v$ by
\begin{equation}
\label{equation:def-uv}
u=\frac{M_1}{N_1},\quad v=\frac{M_2}{N_2},
\end{equation}
where
$$
\begin{array}{lll}
M_1&=&252105 t_5t_7^5+
  11025  s_3^4t_7^4 -172053  s_3^3 t_5^2 t_7^3  -25  s_3^2 t_5 (1080  s_3^5 - 361 t_5^3)t_7^2\\
&&  +
  9504  s_3^6 t_5^3t_7
     -7200  s_3^5 t_5^5, \\
N_1&=&-50421  s_3t_7^5 -180075 t_5^2t_7^4+ 11340  s_3^4 t_5t_7^3+
  129555  s_3^3 t_5^3t_7^2\\
  &&+
     50  s_3^2 t_5^2 (432  s_3^5 + 125 t_5^3)t_7
     -8640  s_3^6 t_5^4,\\
M_2&=&-8823675 
   t_5^2t_7^6 -1065015  s_3^4 t_5t_7^5
     -15  s_3^3 (16875  s_3^5 - 426643 t_5^3)t_7^4\\
     &&+
  25  s_3^2 t_5^2 (45576  s_3^5 - 20825 t_5^3)t_7^3
     -27  s_3 t_5^4 (16664  s_3^5 + 9375 t_5^3) t_7^2\\
     &&+     540000  s_3^5 t_5^6t_7+
  48  s_3^4 t_5^5 (1728  s_3^5 + 3125 t_5^3), \\
N_2&=&1555848  s_3^3 t_5t_7^5+
  525  s_3^2 (675  s_3^5 - 1666 t_5^3)t_7^4 -210  s_3 
   t_5^2 (1548  s_3^5 + 48125 t_5^3)t_7^3\\&&+
     945 t_5^4 (232  s_3^5 - 3125 t_5^3)t_7^2+ 6300000  s_3^4 t_5^6t_7+
  560  s_3^3 t_5^5 (1728  s_3^5 + 3125 t_5^3).\\
  \end{array}
  $$
Actually a direct computation shows that
$$
{\tilde f}_0\left(\frac{M_1}{N_1},\frac{M_2}{N_2},0\right)\equiv
\partial_x{\tilde f}_0\left(\frac{M_1}{N_1},\frac{M_2}{N_2},0\right)\equiv
\partial_y{\tilde f}_0\left(\frac{M_1}{N_1},\frac{M_2}{N_2},0\right)\equiv
\partial_z{\tilde f}_0\left(\frac{M_1}{N_1},\frac{M_2}{N_2},0\right)\equiv0\>({\rm mod}\>\tilde\delta_0),
$$
where ${\tilde f}_0={\tilde f}|_{t_1=t_2=t_3=t_4=0}=x^3y+y^3+t_7x+t_5x^2+s_3y^2-z^2$.
As a consequence,
$(u^{\ast},v^{\ast},0;0,0,0,0,s_3^{\ast},t_5^{\ast},t_7^{\ast})\in{\tilde C}_{\tilde f}$
and therefore $(0,0,0,0,s_3^{\ast},t_5^{\ast},t_7^{\ast})\in  C_{\tilde f}$.
We have thus proved $D_{\tilde f}\subseteq C_{\tilde f}$.
Then (iv) follows.
[]

\begin{definition}
\label{definition:def3-1}
The polynomial ${\tilde\delta}(s_3,\tau)$ is called  the discriminant of $\tilde{\cal F}$.
\end{definition}

\begin{remark}
\label{remark:discriminant}

(i)
Let $\delta_{E_7}=\delta_{E_7}(p_0,p_1,q_0,\cdots,q_4)$ be the discriminant of
the root system of type $E_7$.
Then by the substitution
(\ref{equation:pqst}),
$\delta_{E_7}$ turns out to be a polynomial of $(s_3,\tau)$
which coincides with $t_7^2{\tilde{\delta}}$ up to a non-zero
constant factor.

(ii) The polynomial $\delta_{\rm ST34}(\tau)$ is regarded as the discriminant of the group ${\rm ST34}$.
This will be explained in the next section.

\end{remark}
        
\subsection{Subfamilies of 1-parameter deformations of ${\tilde {\cal F}}$}
\label{subsection:3-3}

We determine the parameter $(s_3,\tau)$ of surfaces corresponding to the parameter $\lambda[i]$
computed in \S2.3.
For this purpose, it is sufficient to
 solve the system of equations
$$
\lambda^{(i)}=(p_0,p_1,q_0,q_1,q_2,q_3,q_4),
$$
where $(p_0,p_1,q_0,q_1,q_2,q_3,q_4)$ are polynomials of $(s_3,\tau)$ defined by
(\ref{equation:pqst}).

The result is given below.

$$
\begin{array}{llll}
{\rm Case }\>E_7{\rm (1)}&&&\\
{\rm Solution}\>1& (s_3,t_1,t_2,t_3,t_4,t_5,t_7)&=&(0,\eta,-3\eta^2,-2\eta^3,0,0,0)\\
{\rm Solution}\>2&  (s_3,t_1,t_2,t_3,t_4,t_5,t_7)&=&(\frac{27}{8}\eta^3,\eta,\frac{15}{4}\eta^2,\frac{65}{8}\eta^3,\frac{27}{16}\eta^4,\frac{675}{32}\eta^5,\frac{2187}{128}\eta^7)\\
&&&\\
{\rm Case}\>E_7{\rm (2)}& (s_3,t_1,t_2,t_3,t_4,t_5,t_7)&=&
 (2\eta^3,\frac{7}{6}\eta,-\frac{10}{3}\eta^2,\frac{52}{27}\eta^3,\frac{8}{3}\eta^4,\frac{8}{9}\eta^5,0).\\
&&&\\
{\rm Case}\>E_7{\rm (3)}&&&\\
{\rm Solution}\>1& (s_3,t_1,t_2,t_3,t_4,t_5,t_7)&=&(\eta^3,\eta,0,2\eta^3,-\eta^4,\eta^5,0)\\
{\rm Solution}\>2& (s_3,t_1,t_2,t_3,t_4,t_5,t_7)&=&(\eta^3,\eta,-3\eta^2,-2\eta^3,2\eta^4,\eta^5,0)\\
{\rm Solution}\>3& (s_3,t_1,t_2,t_3,t_4,t_5,t_7)&=&(-\frac{1}{8}\eta^3,\eta,-\frac{9}{4}\eta^2,-\frac{11}{8}\eta^3,\frac{11}{16}\eta^4,
 \frac{23}{32}\eta^5,-\frac{9}{128}\eta^7)\\
&&&\\
{\rm Case}\>E_7{\rm (4)}&&&\\
{\rm Solution}\>1& (s_3,t_1,t_2,t_3,t_4,t_5,t_7)&=&(4\eta^3,2\eta,-3\eta^2,6\eta^3,-8\eta^4,4\eta^5,0)\\
{\rm Solution}\>2& (s_3,t_1,t_2,t_3,t_4,t_5,t_7)&=&(-2\eta^3,2\eta,-6\eta^2,-4\eta^3,\eta^4,10\eta^5,0)\\
{\rm Solution}\>3& (s_3,t_1,t_2,t_3,t_4,t_5,t_7)&=&(16\eta^3,2\eta,-15\eta^2,-22\eta^3,64\eta^4,
64\eta^5,0)\\
&&&\\
{\rm      Case}\>E_7{\rm (5)}&&&\\
{\rm Solution}\>1& (s_3,t_1,t_2,t_3,t_4,t_5,t_7)&=&(2\eta^3,\frac{5}{2}\eta,-6\eta^2,4\eta^3,-8\eta^4,8\eta^5,0)\\
{\rm Solution}\>2& (s_3,t_1,t_2,t_3,t_4,t_5,t_7)&=&(50\eta^3,\frac{5}{2}\eta,-30\eta^2,-60\eta^3,280\eta^4,392\eta^5,0)\\
{\rm Solution}\>3& (s_3,t_1,t_2,t_3,t_4,t_5,t_7)&=&(-\frac{25}{64}\eta^3,\frac{5}{2}\eta,-\frac{105}{16}\eta^2,\frac{85}{64}\eta^3,
 -\frac{1445}{256}\eta^4,
 \frac{11783}{1024}\eta^5,\frac{30375}{16384}\eta^7)\\
&&&\\
{\rm      Case}\rm E_7{\rm (6)}&&&\\
{\rm Solution}\>1& (s_3,t_1,t_2,t_3,t_4,t_5,t_7)&=&(0,2\eta,-3\eta^2,2\eta^3,0,0,0)\\
{\rm Solution}\>2& (s_3,t_1,t_2,t_3,t_4,t_5,t_7)&=&(54\eta^3,2\eta,-30\eta^2,-52\eta^3,297\eta^4,378\eta^5,0)\\
&&&\\
{\rm      Case }\>E_7{\rm(7)}&&&\\
{\rm Solution}\>1& (s_3,t_1,t_2,t_3,t_4,t_5,t_7)&=&(0,\eta,0,0,0,0,0)\\
{\rm Solution}\>2& (s_3,t_1,t_2,t_3,t_4,t_5,t_7)&=&(\frac{243}{8}\eta^3,\eta,-\frac{81}{4}\eta^2,-\frac{135}{8}\eta^3,\frac{2187}{16}\eta^4,
 \frac{2187}{32}\eta^5,\frac{19683}{128}\eta^7)\\
 \end{array}
 $$

As a consequence, we have the following theorem.

\begin{theorem}
For each corank one subdiagram ${\cal R}_0$  of the root system of type $E_7$,
the family $\tilde{\cal F}$ contains a surface with isolated singular points
which are  simple singularities of the type corresponding
 to ${\cal R}_0$.
\end{theorem}

\begin{remark}
\label{remark:remark-3}
Suggested by the article of Nagano and Shiga \cite{NS},
the author constructed the polynomial  introduced in
(\ref{equation:polynomial-algebraicFrobenius}).
Actually the expression of the polynomial given in \cite{NS} 
naturally leads us to define a family of surfaces of ${\bf C}^3$
 which is a deformation of 
a simple singularity of type $E_6$.
The polynomial $f_{(s_3,\tau)}$  is related with an algebraic potential of type $E_7$ constructed in \cite{Se3}
which  will be given in \S5.
\end{remark}


\section{Basic invariants, the discriminant of ST34 and  a deformation of a simple singularity of type $E_7$}
\label{section:four}

In this section, we l introduce a family ${\cal F}$ consisting of surfaces ${\tilde{f}}_{(0,\tau)}=0$.
Clearly, the surface   ${\tilde{f}}_{(0,\tau)}=0$ of ${\bf C}^3$ has a singular point
if and only if $t_7\cdot \delta_{\rm ST34}(\tau)=0$.
We  show in \S4.2 an identification of $\delta_{\rm ST34}(\tau)$ with the discriminant of the group ST34.
In \S4.3, wel study the adjacency relation of the family ${\cal F}$.

\subsection{Basic invariants of ST34}
\label{subsection:4-1}

 It is known that ST34 is generated by six linear actions on ${\bf C}^6$ with the coordinate system $(x_1,\cdots,x_6)$.
Concretely, ST34 can be generated by  $P_1,P_2,P_3,Q_1,R_1,R_2$ given at p.298 in \cite{ST}.
The   hyperplanes fixed by the generators   are
$$
\begin{array}{ll}P_1 & : 
x_2-x_3=0,\\
P_2&: x_3-x_4=0,\\
P_3&:x_4-x_5=0,\\
Q_1&: x_1-x_2=0,\\
R_1&: x_1-\omega x_2=0,\\
R_2&: x_1+x_2+x_3+x_4+x_5+x_6=0,\\
\end{array}
$$
where $\omega=(-1+\sqrt{3}i)/2$ is a cubic root of unity.
Let $P_4$ be the pseudo-reflection leaving the hyperplane $x_5-x_6=0$ fixed.
Throughout this paper, the complex reflection group $G(3,3,6)$ 
is identified with the subgroup of ST34 generated by
$P_1,P_2,P_3,P_4,Q_1,R_1$.
For a moment, let $S={\bf C}[x_1,\cdots,x_6]$ be the polynomial ring.
As the basic invariants of $G(3,3,6)$, we may take
$p_{3j}\>(j=1,2,3,4,5)$ and $q_6$ defined by
$$
p_{3j}=x_1^{3j}+x_2^{3j}+x_3^{3j}+x_4^{3j}+x_5^{3j}+x_6^{3j}\>(j=1,2,3,4,5),
\>q_6=x_1x_2x_3x_4x_5x_6.
$$

We recall the result of \cite{OS} to continue our arguments.
Let $m_j\>(j=1,2,3,4,5,7)$ be the polynomials of
$p_{3j}\>(1\le j\le 5)$ and $s_6$ defined in \cite{OS}, \S3.
Then it is shown there that  $S^{{\rm ST34}}={\bf C}[m_1,\cdots,m_5,m_7]$.

\subsection{A subfamily of $\tilde{\cal F}$ and the discriminant of ST34}
\label{subsection:4-2}

Putting
$$
f_{\tau}=
{\tilde f}_{(0,\tau)},
$$
we introduce
a surface of ${\bf C}^3(x,y,z)$ by
$$
S(\tau):f_{\tau}(x,y,z)=0,
$$
and the family
$$
{\cal F}=\{S(\tau)|\>\tau\in{\bf C}^6\}.
$$

Similarly to  \S\ref{subsection:3-2},
we define
$$
\begin{array}{l}
{\tilde C}_f=\{(x,y,z;\tau)\in{\bf C}^9|\>f=\partial_xf=\partial_yf=\partial_zf=0\},\\
C_f=\{\tau\in{\bf C}^6|\>\exists (x,y,z)\in{\bf C}^3\>{\rm s.t.}\>(x,y,z;\tau)\in{\tilde C}_f\}.\\
\end{array}
$$

The following lemma is proved by a direct computation.

\begin{lemma}
\label{lemma:st34-dsc}
(i) As a polynomial of $t_7$,
$\delta_{\rm ST34}$ has the form
$$
\delta_{\rm ST34}(\tau)=t_7^6+\sum_{j=1}^6A_j(t_1,t_2,t_3,t_4,t_5)t_7^{6-j}
$$

(ii) $C_f=\{\tau\in{\bf C}^6|\>t_7\cdot \delta_{\rm ST34}(\tau)=0\}.$
\end{lemma}

The following lemma is already shown in \cite{Se3}.
\begin{lemma}
\label{lemma:lemma-3}
If   $t_1,t_2,t_3,t_4,t_5,t_7$ are written as polynomials of
$m_1,\cdots,m_5,m_7$ in the form below,
then $\delta_{\rm ST34}$ coincides with the discriminant of ST34 up to a constant factor.

{\footnotesize

\begin{equation}
\label{equation:ST34-CS}
\left\{
\begin{array}{lll}
t_1 &=& -m_1/3, \\
t_2 &=& \frac{1}{375} (-261 m_1^2 + 136 m_2), \\
     t_3 &=& 2 (4332933 m_1^3 - 3218508 m_1 m_2 + 182200 m_3)/
    35008875, \\
     t_4 &=& 32 (13040967602916 m_1^4 - 19382588618832 m_1^2 m_2 + 
        3749217437791 m_2^2\\
        && + 
               2954928140625 m_1 m_3 - 362524562500 m_4)/
    3025275887390625, \\
     t_5 &=& -\frac{32}{
        4487084376913286386734375} (8914891074037771985316 m_1^5 - 
            15462298032201009500832 m_1^3 m_2\\
            && + 
      5103760684853498949891 m_1 m_2^2 + 
            2363633130501481337325 m_1^2 m_3\\
            && - 
      591120664646165159200 m_2 m_3 - 
            359933922744912562500 m_1 m_4 + 31067730199334950000 m_5), \\
     t_7 &=& 256 (14871910397621845428216716224661671265424 m_1^7 \\
     &&- 
               22594055415192719865974501421984515475072 m_1^5 m_2\\
               && + 
               39247738767940784266429333042398888460122 m_1^3 m_2^2\\
               && - 
               21532675581751266617292332369740068234849 m_1 m_2^3\\
               && - 
               18871924703475618698688320687560320957450 m_1^4 m_3\\
               && - 
               
        7532523232522082240537430915748704681975 m_1^2 m_2 m_3\\
        && + 
               5405897134435116004179853550160817358175 m_2^2 m_3 \\
               &&+ 
               938153790719655723498182223303588515625 m_1 m_3^2\\
               && + 
               9505439156792034835091871022754851500000 m_1^3 m_4\\
               && + 
               3847577219357879613235041793857467250000 m_1 m_2 m_4\\
               && - 
               580393753984805908080293175724529687500 m_3 m_4 \\
               &&- 
               2081725637615256911464889266041895725000 m_1^2 m_5 \\
               &&- 
               679073398830294573402186747103227712500 m_2 m_5 \\
               &&+ 
        55655256504728944788956726765978125000 
                 m_7)/
    2613546665875782334063922888546747071751953125. \\
    \end{array}
    \right.
    \end{equation}
}
\end{lemma}

\vspace{5mm}
{\bf Outline of Proof.}
We explain the outline of the proof of Lemma \ref{lemma:lemma-3}.
To obtain (\ref{equation:ST34-CS}), we are deeply indebted to the software Mathematica.

The weight of $t_j$ is $j/7$ $(j=1,2,3,4,5,7)$ and the homogeneous  degree of $m_{j}$ is $6j$ $(j=1,2,3,4,5,7)$.
Noting these, we may assume that each $t_j$ is a ST34-invariant polynomial of degree $6j$.
Then
\begin{equation}
\label{equation:ttmm}
\left\{
\begin{array}{lll}
t_1&=&c_{11}m_1,\\
t_2&=&c_{21}m_{1}^2+c_{22}m_2,\\
t_3&=&c_{31}m_1^3+c_{32}m_1m_2+c_{33}m_3,\\
t_4&=&c_{41}m_1^4+c_{42}m_1^2m_2+c_{43}m_2^2+c_{44}m_1m_3+c_{45}m_4,\\
t_5&=&c_{51}m_1^5+c_{52}m_1^3m_2+c_{53}m_1m_2^2+c_{54}m_1^2m_3+c_{55}m_2m_3+c_{56}m_1m_4+c_{57}m_5,\\
t_7&=&c_{61}m_1^7+c_{62}m_1^5m_2+c_{63}m_1^3m_2^2+c_{64}m_1m_2^3+c_{65}m_1^4m_3+c_{66}m_1m_3^2\\
&&+c_{67}m_1^2m_2m_3
+m_{68}m_2^2m_3
+c_{69}m_1^3m_4+c_{6,10}m_1m_2m_4+c_{6,11}m_3m_4\\
&&+c_{6,12}m_1^2m_5+c_{6,13}m_2m_5+
c_{6,14}m_7,\\
\end{array}
\right.
\end{equation}
where all the $c_{ij}$ are constants.
We now recall the Saito matrix $M_{34}$ constructed by D. Bessis and J. Michel for the basic invariants
$\mu_6,\mu_{12},\cdots,\mu_{30},\mu_{42}$.  For the details, see \cite{OS}.
As we noted before that $\mu_{6j}$ coincides with $m_j$ up to a constant which is explicitly given 
in \cite{OS}, Corollary 1.
Then $\det(M_{34})$ is regarded as a polynomial of $m_1,m_2,\cdots,m_5,m_7$.
On the other hand, by the substitution (\ref{equation:ttmm}), $\delta_{\rm ST34}$ is also
regarded as a polynomial of  $m_1,m_2,\cdots,m_5,m_7$.
The comparison of the coefficients of both sides of  the equation
\begin{equation}
\label{equation:m34-saitomatrix}
\det(M_{34})=\delta_{\rm ST34}
\end{equation}
implies a system of algebraic equations of $c_{ij}$'s.
Solving this system, we determine these constants and the
 result is shown in (\ref{equation:ST34-CS}).

\vspace{5mm}
If we regard $t_1,t_2,\cdots$ as polynomials of $m_1,m_2,\cdots$ by (\ref{equation:ST34-CS}),
$t_1,t_2,\cdots$ are invariants of $S^{\rm ST34}$.
It is clear from the definition that $t_1,t_2,\cdots$ actually become basic invariants.

\begin{remark}
There are many choices of the basic invariants of $S^{\rm ST34}$.
Historically they are Conway and Sloane \cite{CS} who first constructed the basic invariants denoted by
$\mu_{6j}\>(j=1,\cdots,5,7)$.
By definition, $m_j$ coincides with $\mu_{6j}$ up to a constant factor (cf. \cite{OS}, Corollary 1).
On the other hand, Terao and Enta \cite{TE} (see also \cite{OT})
introduced the basic invariants  $f_j\>(1\le j\le 6)$.
There is another choice of basic invariants.
We now explain its meaning briefly.
Flat structure and flat coordinates for a well-generated complex reflection group were introduced in \cite{KMS2}.
By definition, the flat coordinates are regarded as the generators of the ring of invariant polynomials.
Flat coordinates for ST34 which is  well-generated were constructed in \cite{KMS1} conjecturally and later
were confirmed (see \cite{KMS2},\cite{Se2}).
Therefore the flat coordinates for ST34 given in \cite{KMS1} also form the  basic invariants.
The fifth choice of basic invariants are  $t_1,t_2,t_3,t_4,t_5,t_7$ which appear naturally 
as coefficients of the defining polynomial of a family of hypersurfaces which is a deformation of a simple
singularity of type $E_7$.

\end{remark}

\subsection{Corank one reflection subgroups of ST34  and deformations of simple singularities}
\label{subsection:4-3}

The group ST34 is generated by the pseudo-reflections
$P_1,P_2,P_3,Q_1,R_1,R_2$ as is explained in \S\ref{subsection:4-1}.
We denote them by $\rho_1,\cdots,\rho_6$ (cf. (\ref{equation:generator-ST34}))
and
 linear forms
$$
\begin{array}{lll}
\ell(\rho_1)&=&x_2-x_3,\\
\ell(\rho_2)&=&x_3-x_4,\\
\ell(\rho_3)&=&x_4-x_5,\\
\ell(\rho_4)&=&x_1-x_2,\\
\ell(\rho_5)&=&x_1-\omega x_2,\\
\ell(\rho_6)&=&x_1+x_2+\cdots+x_6
\\
\end{array}
$$
corresponding to the pseudo-reflections $\rho_1,\cdots,\rho_6$.

We explain the purpose of this subsection.
Let ${\tt x}=(x_1,\cdots,x_6)$
and take one of the generators $\rho_j\>(1\le j\le 6)$, say $\rho_i$.
We choose a vector ${\tt x}[i]$
with the condition
$$
\ell(\rho_k)({\tt x})=0\>(1\le k\le 6,\>k\not=i).
$$
Noting that each $t_k$ is regarded as  a polynomial of $(x_1,\cdots,x_6)$
by (\ref{equation:ST34-CS}),
we define $t_k[i]={t_k}|_{{\tt x}\to{\tt x}[i]}$.
Moreover we put
$\tau[i]=(t_1[i],\cdots,t_5[i],t_7[i])$.
Our purpose of this subsection is to determine $\tau[i]$, the polynomial $f_{{\rm ST34},i}=f_{\tau[i]}$
and the singular points of the surface $S(\tau[i])$.
To accomplish this purpose, we need the equations in (\ref{equation:ST34-CS})
and the explicit forms of basic ST34-invariants of  given in \cite{OS}
(they are  firstly done by Conway and Sloane \cite{CS}).
As was remarked at the end of \S\ref{subsection:4-1},
$m_1,\cdots,m_5,m_7$ are polynomials of $p_3,\cdots,p_{15},q_6$,
so are $t_1,\cdots,t_5,t_7$.
As a consequence, the determination of the values of $t_j$'s at some ${\tt x}_-0\in{\bf C}^6$
are performed by filst determining $p_j$'s and $q_6$ at ${\tt x}_0$ adn then
doing $t_j$'s at the values of   $p_j$'s and $q_6$ thus obtained.
The result is given as follows.

\vspace{5mm}
{\bf Case ST34 (1): $\rho_1$}

The reflection group generated by $\rho_j\>(1<j\le 6)$ is $W(A_3)\times W(A_2)$.

In this case,
$$
\begin{array}{rll}
{\tt x}[1]&=&(0,0,m,m,m,-3m),\\
(p_3,p_6,p_9,p_{12},p_{15},q_6)&=&( -24m^3,  732m^6, -19680m^9, 
 531444m^{12}, -14348904m^{15},  0).\\
 \end{array}
 $$
 Then, putting $\eta=36m^6$, we find that
 $$
 \tau[1]=( -14\eta,  -210\eta^2,  -182\eta^3, -5103\eta^4, 
 -57834\eta^5,  -118098\eta^7)
 $$
and $S(\tau[1])$ is defined by
$$
f_{{\rm ST34},1}=-118098 \eta^7 x - 57834 \eta^5 x^2 - 182 \eta^3 x^3 - 14 \eta x^4 + 
 x^3 y + (-5103 \eta^4 x - 210 \eta^2 x^2) y + y^3 - z^2=0.
 $$
 If $\eta\not=0$, this surface in the  $(x_,y_,z)$-space has a simple singularity of type $A_3$ at 
 $(9\eta^2,-144\eta^3,0)$ and has a simple singularity of type $A_2$ at 
 $(-243\eta^2,-2916\eta^3,0)$.

 \vspace{5mm}
 {\bf Case ST34 (2): $\rho_2$}

The reflection group generated by $\rho_j\>(1\le j\le 6,\>j\not=2)$ is $G(3,3,4)\times W(A_1)$.

In this case,
$$
\begin{array}{rll}
{\tt x}[2]&=&(0,0,0,m,m,-2m),\\
(p_3,p_6,p_9,p_{12},p_{15},q_6)&=&( -6m^3,  66m^6, -510m^9, 
 4098m^{12}, -32766m^{15},  0).\\
 \end{array}
 $$
Then, putting $\eta=36m^6$, we find that
 $$
 \tau[2]=( -2\eta,-3\eta^2,-2\eta^3 0,0,0)
 $$
and $S(\tau[2])$ is defined by
$$
f_{{\rm ST34},2}=-2 \eta^3 x^3 - 2 \eta x^4 - 3 \eta^2 x^2 y + x^3 y + y^3 - z^2.
 $$
 Since $f_{{\rm ST34},2}$ coincides with (\ref{equation:alpha-6case}),
 we conclude 
  that
if $\eta\not=0$, this surface in the  $(x,y,z)$-space has a simple singularity of type $D_5$ at the origin
 and 
  has a simple singularity of type $A_1$ at 
 $(-9\eta^2,-18\eta^3,0)$.

  \vspace{5mm}
 {\bf Case ST34 (3): $\rho_3$}

The reflection group generated by $\rho_j\>(1\le j\le 6,\>j\not=3)$ is ST33.
(ST33 is the complex reflection group No.33 in the list of \cite{ST}.)

In this case,
$$
\begin{array}{rll}
{\tt x}[3]&=&(0,0,0,0,m,-m),\\
(p_3,p_6,p_9,p_{12},p_{15},q_6)&=&
( 0,  2m^6,  0,  2m^{12},  0,  0).\\
 \end{array}
 $$
Then, putting $\eta=-4m^6$, we find that
 $$
 \tau[3]=( \eta, 0,0,0,0,0)
 $$
and $S(\tau[3])$ is defined by
$$
f_{{\rm ST34},3}=\eta x^4+x^3y+y^3-z^2=0.
$$
 Since $f_{{\rm ST34},3}$ coincides with (\ref{equation:alpha-7case}),
 we conclude  that
if $\eta\not=0$, this surface in the  $(x_,y_,z)$-space has a simple singularity of type $E_6$ at  the origin.

  \vspace{5mm}
 {\bf Case ST34 (4): $\rho_4$}

The reflection group generated by $\rho_j\>(1\le j\le 6,\>j\not=4)$ is $W(A_5)$.

$$
\begin{array}{rll}
{\tt x}[4]&=&(\omega m, m, m, m, m, -(4 + \omega)m),\\
(p_3,p_6,p_9,p_{12},p_{15},q_6)&=&
( -12 (4 + 3 \omega)m^3,  6  (253 + 420 \omega)m^6,
 -12  (-878 + 8109 \omega)m^9,\\
&&    6  (-676871 + 212520 \omega)m^{12}, 
 12  (21762464 + 10377273 \omega)m^{15}, \\
&&- \omega (4 + \omega)m^6).\\
\end{array}
$$
Then, putting $\eta=2^2\cdot 3^4m^6$,
we find that
$$
\tau[4]= (-(4+5\omega)\eta, 9(5+\omega)\eta^2,
-189\omega\eta^3,
 243(3+2\omega)\eta^4,
  -729(-2+\omega)\eta^5,
   -6561\omega\eta^7).$$
   and $S(\tau[4])$ is defined by
$$
\begin{array}{lll}
f_{{\rm ST34},4}&=&-6561\omega\eta^7 x-729(-2+\omega)\eta^5 x^2-189\eta^3\omega x^3-(4+5\omega)\eta x^4\\
&&+
(243(3+2\omega)\eta^4x+9(5+\omega)\eta^2 x^2)y+y^3-z^2.\\
\end{array}
$$
If $\eta\not=0$, this surface in the  $(x_,y_,z)$-space has a simple singularity of type $A_5$ at  
the point $(-9(1+2\omega)\eta^2,27(2+\omega)\eta^3,0)$.

   \vspace{5mm}
 {\bf Case ST34 (5): $\rho_5$}

The reflection group generated by $\rho_j\>(1\le j\le 6,\>j\not=5)$ is $W(A_4)\times W(A_1)$.

$$
\begin{array}{rll}
{\tt x}[5]&=&( m, m, m, m, m, -5m),\\
(p_3,p_6,p_9,p_{12},p_{15},q_6)
&=&
( -120m^3, 15630m^6,  -1953120m^9, 
 244140630m^{12}, \\
 && -30517578120m^{15},  -5m^6).\\
\end{array}
$$
Then, putting $\eta=2^2\cdot 3^3m^6$,
we find that
$$
\tau[5]=
( -70 \eta,  -10395 \eta^2,  402570 \eta^3, 
 13063680 \eta^4, -838688256 \eta^5,  161243136000 \eta^7)
 $$
    and $S(\tau[5])$ is defined by
$$
\begin{array}{l}
f_{{\rm ST34},5}=
161243136000 \eta ^7 x - 838688256 \eta ^5 x^2 + 402570 \eta ^3 x^3 - 
 70 \eta  x^4 + x^3 y \\
\hspace{30mm} + (13063680 \eta ^4 x - 10395 \eta ^2 x^2) y + y^3 - z^2=0.\\
 \end{array}
 $$
 
 If $\eta\not=0$, this surface in the  $(x_,y_,z)$-space has a simple singularity of type $A_4$ at  
the point $(-1728\eta^2,-139968\eta^3,0)$
 has a simple singularity of type $A_1$ at  
the point $(3375\eta^2,-109350\eta^3,0)$.

 {\bf Case ST34 (6): $\rho_6$}

The reflection group generated by $\rho_j\>(1\le j< 6)$ is $G(3,3,5)$.

$$
\begin{array}{rll}
{\tt x}[6]&=&( 0,0,0,0,0,, m),\\
(p_3,p_6,p_9,p_{12},p_{15},q_6)&=&
( m^3, m^6,  m^9, m^{12},m^{15},0).\\
\end{array}
$$
Then, putting $\eta=m^6/3$,
we find that
$$
\tau[6]=
( -\eta, -3 \eta^2, 2 \eta^3,0,0,0)
$$
    and, by changing the coordinate $(x,y,z)=(x,y_1+\eta x,z)$,
     $S(\tau[6])$ is defined by
$$
f_{{\rm ST34},6}=x^3y_1+3\eta xy_1^2+y_1^3-z^2=0.
$$
 Since $f_{{\rm ST34},6}$ coincides with (\ref{equation:alpha-1case}),
 we conclude that
 if $\eta\not=0$, this surface in the  $(x,y_1,z)$-space has a simple singularity of type $D_6$ at
 the origin (see (\ref{equation:alpha-1case})).

\vspace{5mm}
Summarizing the computation above, we obtain the following theorem.

\begin{theorem}
The family  ${\cal F}_{\rm ST34}$ of surfaces in ${\bf C}^3$ has the following properties.

(i) $S(0,0,0,0,0,0)$ has a unique isolated singular point at the origin whose type is $E_7$.

(ii) The parameters $t_1,t_2,t_3,t_4,t_5,t_7$ are regarded as basic invariants of the group ${\rm ST34}$.

(iii) If $\delta_{\rm ST34}(\tau)$ is the discriminant of ST34,
then 
the surface $S(\tau)$ is smooth if and only if $t_7\cdot \delta_{\rm ST34}\not=0$.

(iv) If the parameter $\tau$ coincides with one of $\tau[i]\>(1\le i\le 6)$,
the corresponding reflection subgroups of ${\rm ST34}$, and the type of simple singularities of $S(\tau[i])$
are summarized in the following table:

$$
\begin{array}{l|l|l}\hline
\tau[1]&W(A_3)\times W(A_2)&A_3+A_2\\\hline
\tau[2]&G(3,3,4)\times W(A_1)&D_5+A_1\\\hline
\tau[3]&{\rm ST33}&E_6\\\hline
\tau[4]&W(A_5)&A_5\\\hline
\tau[5]&W(A_4)\times W(A_1)&A_4+A_1\\\hline
\tau[6]&G(3,3,5)&D_6\\\hline
\end{array}
$$

(v) If $t_7=0$ and $\delta_{{\rm ST34}}\not=0$, then the surface
$S(\tau)$ has a unique simple singular point  at  $(x.y,z)=(0,0,0)$
whose type is $A_1$.

\end{theorem}

{\bf Proof.}
Since (i)-(iv) are proved before the statement of the theorem, it is sufficient to show (v).

To prove (v), we assume that $t_7=0,\>\delta_{\rm ST34}\not=0$.
By a direct computation, we find that
$$
\delta_{\rm ST34}|_{t_7=0}=t_4^3\varphi(t_1,\cdots,t_5).
$$
As a result, if $\delta_{\rm ST34}|_{t_7=0}\not=0$, then $t_4\not=0$.

We first consider the case $t_5\not=0$.
Since
$$
f_{\tau}|_{t_7=0}=y^3+x^3y+t_5x^2+t_3x^3+t_1x^4+y(t_4x+t_2x^2)-z^2,
$$
putting $x_1=x+\frac{t_4}{2t_5}y$, we find that
$$
f_{\tau}|_{t_7=0}=x_1^2\psi_1+y^2\psi_2-z^2,
$$
where $\psi_j=\psi_j(x_1,y)$ is a polynomial of $x_1,y$ such that $\psi_j(0,0)\not=0\>(j=1,2)$.
Then the surface $S(\tau)|_{t_7=0}$ has a simple singularity of type $A_1$ at $(x,y,z)=(0,0,0)$.

We next consider the case $t_5=0$.
In this case, putting $x_1=x+y,y_1=x-y$, we find that
$$
f_{\tau}|_{t_7=0}=x_1^2\psi_3+y_1^2\psi_4-z^2,
$$
where $\psi_j=\psi_j(x_1,y_1)$ is a polynomial of $x_1,y_1$ such that $\psi_j(0,0)\not=0\>(j=3,4)$.
Then the surface $S(\tau)|_{t_7=0}$ has a simple singularity of type $A_1$ at $(x,y,z)=(0,0,0)$.

If $t_5=t_4=0$, then $\delta_{\rm ST34}|_{t_7=0}=0$ and it is not necessary to consider.

The remaining problem is to show that there is no singular point of $S(\tau)$ except the origin.
It is provable to this problem affirmatively. 
Since we need elementary but a little lengthy argument to accomplish it, we omit to prove it. []

\begin{remark}
The concrete expression of the invariants $m_j\>(j=1,2,3,4,5,7)$ by the $G(3,3,6)$-invariants
$p_{k}\>(k=3,6,9,12,15)$ and $s_6$ given in \cite{OS}
help us to
compute the parameter $\tau[j]\>(j=1,2,3,4,5,6)$ in this subsection.

\end{remark}


\section{Algebraic potentials }
\label{section:five}

We  explain the definition of an algebraic potential which is 
needed to show an algebraic potential
related with the family  $\tilde{\cal F}$ of surfaces introduced in \S\ref{section:three}.

\subsection{Definition of algebraic potentials}
\label{subsection:5-1}

We formulate the notion of a potential
which is a solution of  the WDVV equation.
For the details, refer to \cite{KMS2} and the references therein.

  Let $F=F_0+F_1$ be a function of
$(x_1,x_2,\cdots,x_n)$.
Here
$F_0$ is  the polynomial defined by
\begin{equation}
\label{equation:def-F_0}
F_0=\left\{\begin{array}{ll}
\displaystyle{\frac{1}{2}x_1x_n^2+\sum_{j=2}^{n/2}x_jx_{n-j+1}x_n}&(n:\>{\rm even})\\
\displaystyle{\frac{1}{2}x_1x_n^2+\sum_{j=2}^{m-1}x_jx_{2m-j}x_n+\frac{1}{2}x_{m}{}^2x_n}&
(n:\>{\rm odd},\>m=(n+1)/2)\\
\end{array}\right.
\end{equation}
and $F_1$ is a function of $(x_1,x_2,\cdots,x_{n-1})$, independent of $x_n$.

We assume that $F$ is weighted homogeneous.
Namely, there are non-zero constants
$w,w_1,w_2,\cdots,w_n$ such that
$EF=wF$, where
$
\displaystyle{
E=\frac{1}{w_n}\sum_{j=1}^nw_jx_j\partial_{x_j}}
$
is the Euler vector field.
The constants
$w_1,w_2,\cdots,w_n$ are assumed to be rational numbers and
$
0<w_1\le w_2\le\cdots\le w_n=1.
$
We define an $n\times n$ matrix $C$ whose $(i,j)$-entry  denoted by $C_{ij}$ is
$\frac{\partial^2F}{\partial x_i\partial x_{n-j}}$.
We need the partial derivatives of $C$, namely,
\begin{equation}
\label{equation:def-tildB}
\tilde{B}_j=\partial_{x_j}C\quad(j=1,2,\cdots,n).
\end{equation}
It is clear that ${\tilde B}_n=I_n$, the identity matrix.

\begin{definition}\label{def1}
The system of partial differential equations
\begin{equation}
\label{equation:def-WDVVeq}
\left\{\begin{array}{l}
EF=wF\\
{}[\tilde{B}_j,\>\tilde{B}_k]=O\quad(j,k=1,2,\cdots,n)\\
\end{array}
\right.
\end{equation}
for $F$ is called the WDVV equation.
If $F$ is its solution with some additional conditions, then  $F$ is the potential 
of a Frobenius manifold and the coordinate  $(x_1,\cdots,x_n)$ is
called the flat coordinate.
\end{definition}

Associated to the matrix $C$, there is another $n\times n$ matrix $T$ defined by
$$
T=EC,
$$
where $E$ is the Euler vector field defined above.

We now briefly explain the relationship between the totality of polynomial potentials and that of real reflection groups.
Let $W$ be an irreducible  finite reflection group of rank $n$
and let ${\cal R}$ be the set of reflections contained in $W$.
For each reflection $s\in{\cal R}$, there is a linear form $L_s$ such that $s(L_s)=-L_s$.
Then $\prod_{s\in{\cal R}}L_s^2$ is called the discriminant of $W$.
Let $I_1,I_2,\cdots,I_n$ be the homogeneous generators of the ring of invariant polynomials.
Then
$\prod_{s\in{\cal R}}L_s^2$ is expressed  as the polynomial of  $I_1,I_2,\cdots,I_n$.
So we denote by $D(m_1,\cdots,m_n)$ the polynomial
such that $\prod_{s\in{\cal R}}L_s^2=D(I_1,\cdots,I_n)$.
Under this preparation, it is known the following.
There is a polynomial potential $F$ with the flat coordinate $(x_1,\cdots,x_n)$ and $C,T$ are 
$n\times n$ matrices defined above.
Then there is a coordinate transformation $\beta(x_1,\cdots,x_n)=(I_1,\cdots,I_n)$
such that 
$\det(T)$ coincides with $D(\beta(x_1,\cdots,x_n))$ up to a constant factor.
By C. Hertling,  this gives a one to one correspondence between
 the totality of polynomial potentials and that of real reflection groups.

\subsection{A special kind of algebraic potentials}
\label{subsection:5-2}

We focus our attention on a special kind of
potentials called algebraic potentials.

We begin with formulating the construction of a candidate of an algebraic potential $F=F_0+F_1$ in the variables $(x_1,...,x_n)$,
where $F_0$ is defined by (\ref{equation:def-F_0}) and $F_1$ is an algebraic function of $(x_1,\ldots,x_{n-1})$
containing an algebraic function $z$ defined by a quadratic equation  below (\ref{equation:def-z}).
We explain the definition of $F_1$ more precisely.
We first assume the existence of a polynomial ${\tilde F}_1(x_1,\cdots,x_{n-1},\zeta)$
of
$x_1,\cdots,x_{n-1},\zeta$, where $\zeta$ is the variable of weight $w_{n-1}/2$ and is independent of $x_1,\cdots,x_{n-1}$.
We also assume that ${\tilde F}_1$ is weighted homogeneous of total weight same  as that of $F_1$, that is,
$w_1+2w_n=w_1+2$.
Then $F_1$ is obtained from ${\tilde F}_1$ by the substitution $\zeta=z$.
In fact,  $F_1$ is defined as
$$
F_1(x_1,\cdots,x_{n-1},z)={\tilde F}_1(x_1,\cdots,x_{n-1},\zeta)|_{\zeta=z},
$$
where $z$ is an algebraic function of $x_1,\cdots,x_{n-1}$ defined by the equation
\begin{equation}
\label{equation:def-z}
z^2+v(x_1,\cdots,x_{n-2})-x_{n-1}=0.
\end{equation}

\begin{remark}
In this subsection, we only define an algebraic potential containing an algebraic function
which is a solution to a quadratic equation of the form
(\ref{equation:def-z}).
Changing the quadratic equation (\ref{equation:def-z}) by an algebraic equation of higher degree,
we can define  algebraic potentials more generally.
\end{remark}

The matrix $C$ is defined in the case of algebraic potentials similar to the polynomial case.
By using $C$, we define
$$
T=EC,
$$
where $E$ is the Euler vector field defined above.

\subsection{An algebraic potential of type $E_7$}
\label{subsection:5-3}

We introduce an algebraic potential which is related with
the family $\tilde{\cal F}$.
As is explained in \S\ref{subsection:5-1},
there is a one to one correspondence between the totality
of polynomial potentials and that of real reflection groups.
Moreover,  to each  real reflection group  whose type is one of  $A,D,E$,
there associates a versal family of surfaces in ${\bf C}^3$  of the corresponding type.
Then it is a basic question whether there is a family 
 of surfaces in ${\bf C}^3$ which are deformations of 
a simple singularity to each algebraic potential or not.
It is worth trying to attack the question whether there exists an algebraic potential  associated to the family $\tilde{\cal F}$
of surfaces in ${\bf C}^3$
which is introduced in \S3  or not.
Its answer was already  given in \cite{Se3} but we state here the result below for the sake of completeness.

\begin{theorem}
\label{theorem:algebraicpotential}
(i)
Define a function $F=F(x_1,\cdots,x_7)$ by
$$
\begin{array}{ll}
&F\\
=&\frac{{x_7}^2 {x_1}}{2}+\frac{{x_4}^2
   {x_7}}{2}+{x_3} {x_5} {x_7}+{x_2} {x_6} {x_7}\\
&+\frac{3939238656 {x_1}^{15}}{1092455}-\frac{10368 {x_2}
   {x_1}^{13}}{7}+\frac{3456 {x_3} {x_1}^{12}}{7}+\frac{18166464 {x_2}^2
   {x_1}^{11}}{539}-\frac{1728 {x_4} {x_1}^{11}}{7}-\frac{19008}{7} {x_2}
   {x_3} {x_1}^{10}\\
   &-\frac{288 {x_5} {x_1}^{10}}{7}+\frac{216576 {x_2}^3
   {x_1}^9}{7}+\frac{643392 {x_3}^2 {x_1}^9}{49}-\frac{4608}{7} {x_2}
   {x_4} {x_1}^9-\frac{288 {x_6} {x_1}^9}{7}-\frac{1728}{7} {x_2}^2
   {x_3} {x_1}^8\\
   &+\frac{864}{7} {x_3} {x_4} {x_1}^8-\frac{432}{7}
   {x_2} {x_5} {x_1}^8+60264 {x_2}^4 {x_1}^7-\frac{115776}{7} {x_2}
   {x_3}^2 {x_1}^7+\frac{72360 {x_4}^2 {x_1}^7}{343}+3024 {x_2}^2
   {x_4} {x_1}^7\\
   &+\frac{144}{7} {x_3} {x_5} {x_1}^7-\frac{432}{7}
   {x_2} {x_6} {x_1}^7+\frac{42528 {x_3}^3 {x_1}^6}{7}+19152 {x_2}^3
   {x_3} {x_1}^6+\frac{36864}{7} {x_2} {x_3} {x_4} {x_1}^6+576
   {x_2}^2 {x_5} {x_1}^6\\
   &-\frac{72}{7} {x_4} {x_5}
   {x_1}^6+\frac{144}{7} {x_3} {x_6} {x_1}^6+\frac{252252 {x_2}^5
   {x_1}^5}{5}+18216 {x_2}^2 {x_3}^2 {x_1}^5+\frac{2844}{7} {x_2}
   {x_4}^2 {x_1}^5+\frac{942 {x_5}^2 {x_1}^5}{245}\\
   &+4824 {x_2}^3
   {x_4} {x_1}^5-\frac{10944}{7} {x_3}^2 {x_4} {x_1}^5+360 {x_2}
   {x_3} {x_5} {x_1}^5+576 {x_2}^2 {x_6} {x_1}^5-\frac{72}{7}
   {x_4} {x_6} {x_1}^5\\
   &-7008 {x_2} {x_3}^3 {x_1}^4+\frac{36}{7}
   {x_3} {x_4}^2 {x_1}^4+23940 {x_2}^4 {x_3} {x_1}^4+1728
   {x_2}^2 {x_3} {x_4} {x_1}^4+420 {x_2}^3 {x_5}
   {x_1}^4\\
   &+\frac{1632}{7} {x_3}^2 {x_5} {x_1}^4+\frac{456}{7} {x_2}
   {x_4} {x_5} {x_1}^4+360 {x_2} {x_3} {x_6}
   {x_1}^4-\frac{12}{7} {x_5} {x_6} {x_1}^4+38220 {x_2}^6
   {x_1}^3+3032 {x_3}^4 {x_1}^3\\
   &-\frac{102 {x_4}^3 {x_1}^3}{7}+4116
   {x_2}^3 {x_3}^2 {x_1}^3+444 {x_2}^2 {x_4}^2 {x_1}^3-\frac{27}{7}
   {x_2} {x_5}^2 {x_1}^3+\frac{54 {x_6}^2 {x_1}^3}{49}+5082 {x_2}^4
   {x_4} {x_1}^3\\
   &+264 {x_2} {x_3}^2 {x_4} {x_1}^3+156 {x_2}^2
   {x_3} {x_5} {x_1}^3-\frac{276}{7} {x_3} {x_4} {x_5}
   {x_1}^3+420 {x_2}^3 {x_6} {x_1}^3+48 {x_3}^2 {x_6} {x_1}^3\\
   &+24
   {x_2} {x_4} {x_6} {x_1}^3+1988 {x_2}^2 {x_3}^3 {x_1}^2-150
   {x_2} {x_3} {x_4}^2 {x_1}^2+\frac{25}{7} {x_3} {x_5}^2
   {x_1}^2+24402 {x_2}^5 {x_3} {x_1}^2\\
   &-544 {x_3}^3 {x_4}
   {x_1}^2+2436 {x_2}^3 {x_3} {x_4} {x_1}^2+357 {x_2}^4 {x_5}
   {x_1}^2-156 {x_2} {x_3}^2 {x_5} {x_1}^2-\frac{3}{7} {x_4}^2
   {x_5} {x_1}^2\\
   &+48 {x_2}^2 {x_4} {x_5} {x_1}^2+84 {x_2}^2
   {x_3} {x_6} {x_1}^2+\frac{156}{7} {x_3} {x_4} {x_6}
   {x_1}^2-\frac{6}{7} {x_2} {x_5} {x_6} {x_1}^2+5439 {x_2}^7
   {x_1}-812 {x_2} {x_3}^4 {x_1}\\
   &-6 {x_2} {x_4}^3 {x_1}+6468
   {x_2}^4 {x_3}^2 {x_1}+147 {x_2}^3 {x_4}^2 {x_1}+54 {x_3}^2
   {x_4}^2 {x_1}+\frac{5}{2} {x_2}^2 {x_5}^2 {x_1}-\frac{5}{14}
   {x_4} {x_5}^2 {x_1}\\
   &+\frac{9}{7} {x_2} {x_6}^2
   {x_1}+1764 {x_2}^5 {x_4} {x_1}+42
   {x_2}^2 {x_3}^2 {x_4} {x_1}+\frac{160}{3} {x_3}^3 {x_5}
   {x_1}+140 {x_2}^3 {x_3} {x_5} {x_1}\\
   &+20 {x_2} {x_3} {x_4}
   {x_5} {x_1}+294 {x_2}^4 {x_6} {x_1}-12 {x_2} {x_3}^2
   {x_6} {x_1}+\frac{9}{7} {x_4}^2 {x_6} {x_1}+36 {x_2}^2 {x_4}
   {x_6} {x_1}-\frac{10}{7} {x_3} {x_5} {x_6} {x_1}\\
   &+\frac{1484
   {x_3}^5}{15}+\frac{1666 {x_2}^3 {x_3}^3}{3}+{x_3}
   {x_4}^3+\frac{{x_5}^3}{84}+7 {x_2}^2 {x_3} {x_4}^2-\frac{1}{2}
   {x_2} {x_3} {x_5}^2+\frac{{x_3} {x_6}^2}{7}+1029 {x_2}^6
   {x_3}\\
   &+\frac{364}{3} {x_2} {x_3}^3 {x_4}+392 {x_2}^4 {x_3}
   {x_4}+\frac{49 {x_2}^5 {x_5}}{5}+21 {x_2}^2 {x_3}^2
   {x_5}-\frac{1}{2} {x_2} {x_4}^2 {x_5}+7 {x_2}^3 {x_4} {x_5}-6
   {x_3}^2 {x_4} {x_5}\\
   &-\frac{20 {x_3}^3 {x_6}}{3}+84 {x_2}^3
   {x_3} {x_6}+8 {x_2} {x_3} {x_4} {x_6}+2 {x_2}^2 {x_5}
   {x_6}+\frac{{x_4} {x_5} {x_6}}{7}+\frac{8}{105}z^5,\\
   \end{array}
$$
where $z$ is an algebraic function of $x_1,x_2,x_3,x_4,x_5,x_6$ defined by
the equation
\begin{equation}
\label{equation:def-ans-z-e7e}
v-{x_6}+z^2=0
\end{equation}
and $v$ is the polynomial of $x_1,\cdots,x_6$  defined by
\begin{equation}
\label{equation:def-ans-z-e7v}
\begin{array}{lll}
v=-36 {x_1}^6-36 {x_1}^4 {x_2}+12 {x_1}^3 {x_3}-231 {x_1}^2
   {x_2}^2-6 {x_1}^2 {x_4}\\
\hspace{10mm}   -84 {x_1} {x_2} {x_3}
  -{x_1}
   {x_5}-49 {x_2}^3-7 {x_2} {x_4}-7 {x_3}^2.\\
   \end{array}
 \end{equation}

In this case,
$w_j=j/7\>(j=1,2,\cdots,7)$
and $F$ is weighted homogeneous and a solution of the WDVV equation.

(ii)
In spite that $\det(T)$ is an algebraic function of $x_1,\cdots,x_7$,
$\det(T)$ is regarded as a polynomial of $x_1,\cdots,x_5,x_7,z$
by eliminating $x_6$  using  the relation (\ref{equation:def-ans-z-e7e}).
Then
by  the coordinate transformation      
$$
  \begin{array}{lll}
{x_1}&=& -3 {c_7}^2 {t_1},\\
{x_2}&=& -\frac{3}{7} {c_7}^4 \left(39
   {t_1}^2+28 {t_2}\right),\\
   {x_3}&=& \frac{3}{7} {c_7}^6 \left(-392
   {s_3}-1581 {t_1}^3-1428 {t_1} {t_2}+1176 {t_3}\right),\\
   {x_4}&=&
   \frac{1}{64} {c_7} \left(-32 {s_3} {t_1}-105 {t_1}^4-120 {t_1}^2
   {t_2}+96 {t_1} {t_3}-16 {t_2}^2+64 {t_4}\right),\\
   {x_5}&=& -\frac{3
   {c_7}^3}{10976} \left(215600 {s_3} {t_1}^2-21952 {s_3} {t_2}+585915
   {t_1}^5+816200 {t_1}^3 {t_2}-646800 {t_1}^2 {t_3}+215600 {t_1} {t_2}^2\right.\\
 &&\left.  -241472 {t_1} {t_4}-241472 {t_2} {t_3}+307328
   {t_5}\right),\\
   {x_7}&=& \frac{1}{843308032}
   (32269440 {s_3}^2 {t_1}-747369560
   {s_3} {t_1}^4-481736640 {s_3} {t_1}^2 {t_2}+408746240 {s_3}
   {t_1} {t_3}\\
   &&+32269440 {s_3} {t_2}^2-180708864 {s_3} {t_4}-1498061973
   {t_1}^7-2850681036 {t_1}^5 {t_2}+2242108680 {t_1}^4 {t_3}\\
   &&-1494739120
   {t_1}^3 {t_2}^2+686109760 {t_1}^3 {t_4}+2058329280 {t_1}^2 {t_2}
   {t_3}-613119360 {t_1}^2 {t_5}-160578880 {t_1} {t_2}^3\\
   &&+408746240
   {t_1} {t_2} {t_4}-613119360 {t_1} {t_3}^2+204373120 {t_2}^2
   {t_3}-301181440 {t_2} {t_5}-301181440 {t_3} {t_4}\\
   &&+843308032{t_7}),\\
   {z}&=& -2352 {c_7}^6 {s_3}\\
   \end{array}
$$
where $c_7$ is a constant satisfying
${c_7}^7=\frac{1}{14112}$,
 $\det(T)$ coincides with the discriminant $\tilde{\delta}$ of the family $\tilde{\cal F}$
up to a non-zero constant factor.

\end{theorem}

{\bf Outline of Proof:}
The proof of the theorem is accomplished by a direct computation.
The method of finding the algebraic potential of the theorem is explained in \cite{Se3} (see also \cite{DS}, \S3).

 Since the function $F$ is given concretely, it is possible to compute $C$
and next ${\tilde B}_j\>(j=1,\cdots,n)$
(cf. (\ref{equation:def-tildB})) starting from $F$.
Then it is an easy job to show $[{\tilde B}_j,\>{\tilde B}_k]=O\>(j,k=1,\cdots,n)$.
In this manner (i) follows.
The proof of (ii) is also accomplished by a direct computation.

\begin{remark}
(i) The author was informed by T. Pouvropoulos for the possibility of the  existence of an algebraic potential
of seven variables
with the weights $w_j=j/7\>(j=1,2,\cdots,7)$.
It is  a conjecture that the potential $F$ in the theorem corresponds to a certain conjugacy class
of $W(E_7)$.

(ii) Results analogous to Theorem \ref{theorem:algebraicpotential}
 above for the cases $E_6$, $E_8$  were treated in \cite{Se3}. 
(See also  Remark \ref{remark:remark-3}.)
In the case $E_6$, the group ST33 plays the role analogous to ST34 in the case $E_7$.
But in the case $E_8$, there is no complex reflection group which plays the role analogus to ST34.

\end{remark}

\noindent
{\bf Acknowledgement}

The author thanks E. Bannai for the discussion on the relation between ST34 and $W(E_7)$,
J. Matsuzawa (Nara Women's Univ) for pointing out the part of the reference \cite{Sd1}
concerning the topic in \S2.3 of this paper
and M. Oura (Kanazawa Univ.) for assisting the research of references and discussions on the group ST34.
This work was partially supported by JSPS KAKENHI Grant Number 17K05269.


\end{document}